\documentclass[a4paper,10pt,twoside]{article}
\usepackage[hmarginratio=1:1]{geometry}
\usepackage[utf8]{inputenc}

\usepackage{amsmath,amsfonts,amssymb,comment,bm}
\usepackage{amsthm}
\usepackage{mathrsfs}
\usepackage{subfigure}
\usepackage{listings}
\usepackage{dsfont}
\usepackage{lscape}
\usepackage{rotating}

\usepackage{tikz}
\usetikzlibrary{arrows,calc,through,backgrounds,matrix,decorations.pathmorphing,positioning}

\usepackage[hidelinks]{hyperref}

\usepackage{lineno}
\DeclareMathAlphabet{\mathpzc}{OT1}{pzc}{m}{it}
\usepackage[mathcal]{euscript}
\usepackage{calligra}
\DeclareMathAlphabet{\mathcalligra}{T1}{calligra}{m}{n}

\usepackage{epstopdf}
\usepackage{graphicx}

\usepackage{titlesec}
\usepackage{titletoc}
\usepackage{appendix}

\makeatletter
\newcommand{\pushright}[1]{\ifmeasuring@#1\else\omit\hfill$\displaystyle#1$\fi\ignorespaces}
\newcommand{\pushleft}[1]{\ifmeasuring@#1\else\omit$\displaystyle#1$\hfill\fi\ignorespaces}
\makeatother

\renewcommand{\thesection}{\arabic{section}}%

\usepackage{etoolbox}
\makeatletter
\patchcmd{\chapter}{\if@openright\cleardoublepage\else\clearpage\fi}{}{}{}
\makeatother

\numberwithin{equation}{section}

\usepackage{fancyhdr}
\newcommand{\np}[1]{
\newpage
\thispagestyle{empty}
\mbox{}
}

\usepackage[style=alphabetic,
            backend=bibtex,
            isbn=false,
            doi=false,
            url=false,
            hyperref=true, 
            firstinits=true,
            block=none]{biblatex}
\renewbibmacro*{journal+issuetitle}{%
  \usebibmacro{journal}%
  \setunit*{\addspace}%
  \iffieldundef{series}
    {}
    {\newunit
     \printfield{series}%
     \setunit{\addspace}}%
  \printfield{volume}%
  \setunit{\addspace}%
  \usebibmacro{issue+date}%
  \setunit{\addcolon\space}%
  \usebibmacro{issue}%
  \setunit{\addcomma\space}%
  \printfield{number}%
  \newunit}

\renewbibmacro*{volume+number+eid}{%
\printfield{volume}
\printfield{number}
}

\renewbibmacro*{date}{\setunit{\addspace}\printdate}

\DeclareFieldFormat*{title}{\mkbibemph{#1}}
\DeclareFieldFormat*{journaltitle}{#1} 
\DeclareFieldFormat*{volume}{vol. \bf #1}
\DeclareFieldFormat[article]{date}{#1}
\DeclareFieldFormat[book,inbook,incollection,thesis,unpublished]{date}{(#1)}
\DeclareFieldFormat*{number}{no. #1}
\DeclareFieldFormat*{pages}{#1}
\renewbibmacro{in:}{}

\newcounter{entryc}[section]
\newcommand{\thetheorem}{\thesection.\arabic{entryc}} 

\makeatletter
\newenvironment{mydef}{%
  \par\addvspace{\baselineskip}\refstepcounter{entryc}\protected@edef\@currentlabel{\thetheorem}%
  \noindent\textbf{Definition \thetheorem.}%
}{\par\addvspace{\baselineskip}}
\makeatother

\makeatletter
\newenvironment{mylem}{%
  \par\addvspace{\baselineskip}\refstepcounter{entryc}\protected@edef\@currentlabel{\thetheorem}%
  \noindent\textbf{Lemma \thetheorem.}%
}{\par\addvspace{\baselineskip}}
\makeatletter
\newenvironment{myprop}{%
  \par\addvspace{\baselineskip}\refstepcounter{entryc}\protected@edef\@currentlabel{\thetheorem}%
  \noindent\textbf{Proposition \thetheorem.}%
}{\par\addvspace{\baselineskip}}
\makeatother
\makeatletter
\newenvironment{mythm}{%
  \par\addvspace{\baselineskip}\refstepcounter{entryc}\protected@edef\@currentlabel{\thetheorem}%
  \noindent\textbf{Theorem \thetheorem.}%
}{\par\addvspace{\baselineskip}}
\makeatother


\makeatletter

\makeatletter
  {\par\addvspace{\baselineskip}}
\makeatother
\makeatletter
\newenvironment{mynot}{%
  \par\addvspace{\baselineskip}\refstepcounter{entryc}\protected@edef\@currentlabel{\thetheorem}%
  \noindent\textbf{Notation \thetheorem.}%
  }
  {\par\addvspace{\baselineskip}}
\makeatother
\makeatletter
  {\par\addvspace{\baselineskip}}
\makeatother
\makeatletter
\newenvironment{mypf}{%
  \noindent\emph{Proof}:%
}{\par\addvspace{\baselineskip}}
\makeatother



\newcommand\blfootnote[1]{%
  \begingroup
  \renewcommand\thefootnote{}\footnote{#1}%
  \addtocounter{footnote}{-1}%
  \endgroup
}

\DeclareMathOperator*{\fC}{\mathscr{C}}

\DeclareMathOperator{\Der}{Der}
\DeclareMathOperator{\Hom}{Hom}

\DeclareMathOperator{\IM}{Im}
\DeclareMathOperator{\Ker}{Ker}

\DeclareMathOperator{\HH}{HH}

\DeclareMathOperator{\Ab}{\mathbf{Ab}}

\DeclareMathOperator{\Set}{\mathbf{Set}}
\DeclareMathOperator{\Top}{\mathbf{Top}}

\DeclareMathOperator{\Cat}{\mathbf{Cat}}

\DeclareMathOperator{\ot}{\wedge}
\allowdisplaybreaks

\title{Hochschild cohomology of ring objects in monoidal categories}
\date{}
\author{Magnus Hellstrøm-Finnsen}

\begin{document}
\maketitle

\begin{abstract}
We define the Hochschild complex and cohomology of a ring object in a monoidal category enriched over abelian groups. We interpret the cohomology groups and prove that the cohomology ring is graded-commutative. 
\end{abstract}


\tableofcontents

\section{Introduction}
With inspiration from the classical definition of Hochschild cohomology (which for example can be found in \cite{hap-89}) we define the similar notion of Hochschild cohomology for ring objects in monoidal categories enriched over abelian groups. Hochschild cohomology was initially studied by Hochschild in \cite{hoc-45} and \cite{hoc-46}. Monoidal categories were introduced as a generalisation (or formalisation) of tensor products, and they have played a fundamental role in the development of category theory ever since. We restate the Hochschild cochain complex in this setting, look at some immediate consequences for lower dimensions and prove that the cohomology ring is graded-commutative. Related work was done in \cite{ams-07}, where the authors studied Hochschild cohomology of abelian monoidal categories by construction a bar complex. 

\section{Preliminaries} 
In this section we discuss monoidal categories, ring objects and their basic properties. Throughout this paper, in a category $\fC$ we denote an object $X$ by $X\in\fC$, and its identity by $1_X$, or only by $1$ if the object is assumed to be known. Morphisms are denoted by arrows $X \to Y$ and natural transformations with arrows such as $F \rightharpoonup G$. 

Let $\Ab$ denote the category of abelian groups and group homomorphisms. A category $\fC$ is said to be an \emph{$\Ab$-enriched category} if the hom-objects in $\fC$ are abelian groups and the composition is bilinear over the integers $\mathbb{Z}$. Note that we do not assume that an $\Ab$-enriched category has a zero object. But if an $\Ab$-enriched category happens to have an initial object, this object is also terminal, hence a zero object, because the zero morphism is in every hom-object. Similarly every finite coproduct in an $\Ab$-enriched category is also a finite product, which often is referred to as a biproduct. If an $\Ab$-enriched category $\fC$ happens to admit finite biproducts, we say that $\fC$ is \emph{additive}. 

\subsection{Monoidal categories} 
We recall the definition of a monoidal category from \cite{mac-98}, 

\begin{mydef} A category $\fC$ is a \emph{monoidal category} if it is equipped with a bifunctor $\ot : \fC \times \fC \to \fC$ (often referred to as the \emph{tensor product}) and an object $I$ in $\fC$ called the \emph{tensor unit}, together with the following natural isomorphisms:  
\begin{itemize}
\item The \emph{associator} $\alpha:( ? \ot ? ) \ot ? \rightharpoonup ? \ot ( ? \ot ? )$ which has components 
\begin{align} 
\alpha_{X,Y,Z}:(X\ot Y)\ot Z \to X\ot (Y\ot Z)
\end{align}
for all objects $X$, $Y$ and $Z$ in $\fC$.
\item The \emph{left unitor} $\lambda: I \ot ? \rightharpoonup ? $  which has components 
\begin{align}
\lambda_X: I \ot X \to X 
\end{align}
for every object $X$ in $\fC$. 
\item The \emph{right unitor} $\rho: ? \ot 1 \rightharpoonup ? $ which has components 
\begin{align}
\rho_X: X \ot I \to X
\end{align}
for every object $X$ in $\fC$.
\end{itemize}
These data should make the \emph{pentagon diagram}
\begin{center}
\begin{tikzpicture}
\matrix(m)[matrix of math nodes,row sep=5em,column sep=-2em,text height=1.5ex,text depth=0.25ex]
{                                  &                                  &(W\ot X) \ot (Y \ot Z)&                                  &                                  \\
 ((W\ot X) \ot Y) \ot Z&                                  &                                  &                                  &W\ot (X \ot (Y \ot Z))\\
                                   &(W\ot (X \ot Y)) \ot Z&                                  &W\ot ((X \ot Y) \ot Z)&                                  \\};
\draw[ ->,font=\scriptsize](m-2-1) edge         node[above,sloped]{$\alpha_{W \ot X,Y,Z}          $} (m-1-3);
\draw[ ->,font=\scriptsize](m-1-3) edge         node[above,sloped]{$\alpha_{W,X,Y \ot Z}          $} (m-2-5);

\draw[ ->,font=\scriptsize](m-2-1) edge         node[left ]{$\alpha_{W,X,Y}\ot 1_Z         $} (m-3-2);
\draw[ ->,font=\scriptsize](m-3-2) edge         node[below]{$\alpha_{W,X\ot Y,Z}           $} (m-3-4);
\draw[ ->,font=\scriptsize](m-3-4) edge         node[right]{$1_W\ot\alpha_{X,Y,Z}          $} (m-2-5);
\end{tikzpicture}
\end{center}
where $W$, $X$, $Y$ and $Z$ are arbitrary objects in $\fC$, and the \emph{triangle diagram}
\begin{center}
\begin{tikzpicture}
\matrix(m)[matrix of math nodes,row sep=3em,column sep=2em,text height=1.5ex,text depth=0.25ex]
{ (X \ot I) \ot Y  &             & X \ot (I \ot Y) \\
                           & X \ot Y &                       \\};
\draw[ ->,font=\scriptsize](m-1-1) edge         node[above]{$ \alpha_{X,I,Y}        $} (m-1-3);                       
\draw[ ->,font=\scriptsize](m-1-1) edge         node[left ]{$ \rho_X \ot 1_Y    $} (m-2-2);                       
\draw[ ->,font=\scriptsize](m-1-3) edge         node[right]{$ 1_X \ot \lambda_Y $} (m-2-2);
\end{tikzpicture} 
\end{center}
where $X$ and $Y$ are arbitrary objects in $\fC$ commutative. We denote the data for this category by $(\fC,\ot,I,\alpha,\lambda,\rho)$. \end{mydef}

Among many examples of monoidal categories we mention the following. The category $\Set$ (of sets and functions) with tensor product given by the cartesian product and the tensor unit given by the one point set $*$ is monoidal. This is even an example of a \emph{cartesian monoidal category}, which is a monoidal category where the monoidal structure is given by the cartesian product and the terminal object is the tensor unit. Another example of a cartesian monoidal category is $\Cat$ (of (small) categories and functors) with products of categories and the terminal category. Similarly any category with coproducts is a monoidal category, where the initial object is tensor unit. Further the category $\Ab$ (of abelian groups and group homomorphisms) with the usual tensor product $\ot=\otimes_{\mathbb{Z}}$ and tensor unit $\mathbb{Z}$. In fact $\Ab$ is even a symmetric monoidal category (defined in Section \ref{sec:centre}). The category $\mathbf{vec}(k)$ (of finite dimensional vector spaces over a field $k$ and linear transformations) is a monoidal category with tensor product $\otimes_k$ and tensor unit $k$. The category $\Top_*$ (of pointed topological spaces and continuous functions) is a monoidal category with the smash product $\wedge$ and the base point as tensor unit. 

\subsection{Ring objects}
Here we will use the term ring object for what many will refer to as monoids (i.e.\ \cite{mac-98}). This terminology is also used in \cite{hps-97}. 

\begin{mydef}\label{def:rngobj} 
Let $(\fC,\ot,I,\alpha,\lambda,\rho)$ be a monoidal category. A \emph{ring object} $R$ is an object in $\fC$ equipped with a \emph{multiplication rule} $\mu_R:R \ot R \to R$ and a \emph{multiplicative unit} $e_R:I \to R$. These morphisms satisfy the following relations: 
\begin{itemize}
\item The \emph{associative relation}: the multiplication rule is associative in the sense that the following diagram commutes 
\begin{center}
\begin{tikzpicture}
\matrix(m)[matrix of math nodes,row sep=2.6em,column sep=2.8em,text height=1.5ex,text depth=0.25ex]
{(R \ot R)\ot R  &               & R \ot(R \ot R)  \\
  R           \ot R  &               & R \ot R             \\
                         &  R            &                         \\};
\draw[ ->,font=\scriptsize](m-1-1) edge         node[above]{$\alpha_{R,R,R}          $} (m-1-3);
\draw[ ->,font=\scriptsize](m-1-1) edge         node[left ]{$\mu_R \ot 1_R           $} (m-2-1);
\draw[ ->,font=\scriptsize](m-2-1) edge         node[below]{$\mu_R                   $} (m-3-2);

\draw[ ->,font=\scriptsize](m-1-3) edge         node[right]{$1_R   \ot \mu_R         $} (m-2-3);
\draw[ ->,font=\scriptsize](m-2-3) edge         node[below]{$\mu_R                   $} (m-3-2);
\end{tikzpicture}
\end{center}
\item The \emph{unitary relation}: the multiplication admits a left unit and a right unit in the sense that the following diagram commutes 
\begin{center}
\begin{tikzpicture}
\matrix(m)[matrix of math nodes,row sep=2.6em,column sep=2.8em,text height=1.5ex,text depth=0.25ex]
{ R                      & I \ot R   & R \ot R        & R \ot I  & R     \\
                         &               & R                  &                      \\
};
\draw[ ->,font=\scriptsize](m-1-1) edge         node[above]{$\lambda^{-1}_R          $} (m-1-2);
\draw[ ->,font=\scriptsize](m-1-2) edge         node[above]{$e_R\ot 1_R              $} (m-1-3);
\draw[ ->,font=\scriptsize](m-1-3) edge         node[right]{$\mu_R                   $} (m-2-3);
\draw[ ->,font=\scriptsize](m-1-1) edge         node[below]{$1_R                     $} (m-2-3);

\draw[ ->,font=\scriptsize](m-1-5) edge         node[above]{$\rho^{-1}_R             $} (m-1-4);
\draw[ ->,font=\scriptsize](m-1-4) edge         node[above]{$1_R\ot e_R              $} (m-1-3);
\draw[ ->,font=\scriptsize](m-1-5) edge         node[below]{$1_R                     $} (m-2-3);
\end{tikzpicture}
\end{center}
\end{itemize}
We denote a ring object as a triple $(R,\mu_R,e_R)$, and often the subscripts are skipped. 
\end{mydef}

In the monoidal category $\Ab$ the ring objects are simply ordinary rings. The ring objects in the monoidal category $\mathbf{vec}(k)$ are (finite dimensional) $k$-algebras. 

\begin{mydef}\label{def:mro}
\sloppy Let $(R,\mu_R,e_R)$ and $(S,\mu_S,e_S)$ be ring objects in a monoidal category $(\fC,\ot,I,\alpha,\lambda,\rho)$. A \emph{morphism of ring objects} $f:R \to S$ is a morphism $f$ in $\fC$ such that 
\begin{align}
&f \circ \mu_R =\mu_S \circ (f \ot f):R \ot R \to S \qquad \text{and} \qquad f \circ e_R = e_S , 
\end{align}
which means that the following two diagrams commute 
\newline
\begin{minipage}{.5\textwidth}
\begin{center}
\begin{tikzpicture}
\matrix(m)[matrix of math nodes,row sep=2.6em,column sep=2.8em,text height=1.5ex,text depth=0.25ex]
{ 
  R \ot R   &  S \ot S    \\
  R         &  S          \\
};
\draw[ ->,font=\scriptsize](m-1-1) edge         node[above]{$f \ot f         $} (m-1-2);
\draw[ ->,font=\scriptsize](m-2-1) edge         node[above]{$f               $} (m-2-2);
\draw[ ->,font=\scriptsize](m-1-1) edge         node[right]{$\mu_R           $} (m-2-1);
\draw[ ->,font=\scriptsize](m-1-2) edge         node[left ]{$\mu_S           $} (m-2-2);
\end{tikzpicture}
\end{center}
\end{minipage}
\begin{minipage}{.5\textwidth}
\begin{center}
\begin{tikzpicture}
\matrix(m)[matrix of math nodes,row sep=.1em,column sep=2.8em,text height=1.5ex,text depth=0.25ex]
{ 
           &  R          \\
       I   &             \\
           &  S          \\
};
\draw[ ->,font=\scriptsize](m-2-1) edge         node[above]{$e_R             $} (m-1-2);
\draw[ ->,font=\scriptsize](m-2-1) edge         node[below]{$e_S             $} (m-3-2);
\draw[ ->,font=\scriptsize](m-1-2) edge         node[right]{$f               $} (m-3-2);
\end{tikzpicture}
\end{center} 
\end{minipage}
\end{mydef}

\subsection{Coherence in monoidal categories}
Throughout this section let $(\fC,\ot,I,\alpha,\lambda,\rho)$ be a monoidal category. We say that $\fC$ is \emph{strict} if $\alpha$, $\lambda$ and $\rho$ all are identities. When confusion may occur, we use the term \emph{weak} for non-strict monoidal categories. Nevertheless, examples like $\Ab$, even $\Set$, which we think of as natural examples of monoidal categories, are not strict. 

However, in weak monoidal categories, as discussed in \cite[Section VII.2]{mac-98} (and originally suggested in \cite{mac-63} and \cite{kel-64}), any formal diagram built up from instances of $\alpha$, $\lambda$ and $\rho$ by $\ot$ commutes. This result will often be referred to as coherence in monoidal categories or the coherence theorem for monoidal categories. The coherence theorem implies that the constructions we can make of objects in $\fC$ by ``moving parentheses with associators'', or ``tensoring with the tensor identity on the left or the right'' are not only isomorphic objects, but all different ways to ``produce'' these constructions are ``equivalent''. ``Equivalent'' in the sense that every diagram containing different procedures of constructing objects using instances of $\alpha$, $\lambda$ and $\rho$ by $\ot$ commutes. 

As a consequence of the coherence theorem it is sometimes usual to not differ between such isomorphic objects, and hence thinking about any monoidal category as strict. However we choose to be loyal to the philosophy that every arrow should start in a distinct object and end in a distinct object, hence we differ between all such coherent cases. 

Later we are going to study long ``chains'' of ring objects tensored together. Hence we introduce the following notation. 

\begin{mynot}
Let $(R,\mu,e)$ be a ring object in $\fC$. We denote 
\begin{align}
R^{\ot k}=(\cdots((R \ot R) \ot R ) \ot R \cdots ) \ot R, 
\end{align}
where $R$ occurs $k$ times and the parentheses are as above (i.e.\ all the left parentheses are grouped together). We use the convention that $R^{\ot 0}$ is the empty symbol. By the coherence theorem all such procedures are equivalent, in the sense that all formal diagrams involving these associators possible tensored with identities commute. We denote this described procedure by 
\begin{align}
&\alpha^{i,j}_k: R^{\ot k} \to ( R^{\ot i} \ot R^{\ot j} ) \ot R^{\ot (k-i-j)} 
\end{align} 
for $0 \leq i \leq k$, $0 \leq j \leq k$ and $0 \leq i+j \leq k$. We remark that $\alpha^{i,j}_k$ in many cases might be the identity, e.g.\ when $k=0,1,2$, but also in many other cases. Now consider again chains of the same ring object tensored together. By 
\begin{align}
&\mu^i_k=(1_{R^{\ot i}} \ot \mu) \ot 1_{R^{\ot(k-i-2)}}:( R^{\ot i} \ot R^{\ot 2} ) \ot R^{\ot (k-i-2)} \to ( R^{\ot i} \ot R ) \ot R^{\ot (k-i-2)} 
\end{align}
for $0 \leq i \leq k-2 $ we denote the multiplication of two objects occurring after a chain of $i$ objects. 
\end{mynot}

\section{The Hochschild complex} 
In this section, we for ring objects and bimodules define the Hochschild complex.

\subsection{The Hochschild cochain complex of a ring object} 
We restate the classical definition of the Hochschild cochain complex from \cite{hoc-45} and \cite{hoc-46} in the setting of ring objects in $\Ab$-enriched monoidal categories.  

\begin{mydef}\label{def:hcc}
Let $(\fC,\ot,I,\alpha,\lambda,\rho)$ be an $\Ab$-enriched monoidal category and let $(R,\mu,e)$ be a ring object in $\fC$. The \emph{Hochschild cochain complex} $C^{\bullet}(R)=(C^k(R),d^k)_{k\in\mathbb{Z}}$ is defined to be the sequence 
\begin{align*}
\cdots \to 0 \xrightarrow{d^{-1}} C^0(R) \xrightarrow{d^{0}} C^1(R) \xrightarrow{d^{1}} C^2(R) \xrightarrow{d^{2}} \cdots 
\end{align*}
that has objects 
\begin{align}
C^k=
\begin{cases}
0                           &\text{for } k<0      \\
\Hom_{\fC}(I,R)             &\text{for } k=0      \\ 
\Hom_{\fC}(R^{\ot k},R) &\text{for } k\geq 1. 
\end{cases}
\end{align}
The differentials $d^k: C^k(R) \to C^{k+1}(R)$ are defined as:  
\begin{itemize}
\item $d^k=0$ for $k<0$. 
\item For $f\in C^0(R)=\Hom_{\fC}(I,R)$ the differential $d^0:\Hom_{\fC}(I,R)\to\Hom_{\fC}(R,R)$ is defined to be 
\begin{align}
d^0(f)=&\mu \circ (1_R\ot f) \circ \lambda^{-1}_R - \mu \circ (f\ot 1_R) \circ \rho^{-1}_R. 
\end{align}
\item For $k \geq 1$ and $f\in C^k(R) = \Hom_{\fC}(R^{\ot k},R)$ the differentials $d^k:\Hom_{\fC}(R^{\ot k},R)\to\Hom_{\fC}(R^{\ot(k+1)},R)$ are defined to be 
\begin{align}
d^k(f)={ {\mu} \circ ({1_R \ot f}) \circ \alpha^{0,1}_{k+1} }  
      &+\sum^{k}_{i=1}(-1)^i[ {f} \circ (\alpha^{i-1,1}_k)^{-1} \circ {\mu^{i-1}_{k+1}} \circ \alpha^{i-1,2}_{k+1} ] \nonumber\\
      &+(-1)^{k+1}[{\mu} \circ {(f \ot 1_R)}]\nonumber
\end{align}\end{itemize}\end{mydef}

Next we prove that this sequence is a cochain complex. The proof spans several pages. 

\begin{mythm}\label{thm:hcc}
Let $(\fC,\ot,I,\alpha,\lambda,\rho)$ be an $\Ab$-enriched  monoidal category and let $(R,\mu,e)$ be a ring object in $\fC$. The sequence $C^{\bullet}(R)=(C^k(R),d^k)_{k\in\mathbb{Z}}$ is a cochain complex, i.e.\ $d^{k+1} \circ d^{k}=0$ for all $k\in\mathbb{Z}$.\end{mythm}

\begin{mypf}
This proof mainly consists of identifying terms in sums with one another, in such a way that they cancel each other. Braces are used to label the different terms. The sign of the term is always considered as a part of the term. We divide the proof into separate cases. 
\\\\
\textbf{The case $k<0$.} This case is obvious since all $d^k=0$. 
\\\\
\textbf{The case $k=0$.} For this case we want to prove if $d^1 \circ d^0 = 0$. Recall $d^1(?)=\mu \circ (1_R\ot ?) +(-1)?\mu + \mu \circ (?\ot 1_R)$. For $f\in\Hom_{\fC}(I,R)$ we get 
\begin{align}
(d^1 \circ d^0)(f)=&\mu \circ (1_R \ot d^0(f)) - d^0(f) \circ \mu + \mu \circ (d^0(f) \ot 1_R)                                                      \nonumber \\
                =&\overbrace{\mu \circ [1_R \ot \mu \circ (f \ot 1_R) \circ \lambda^{-1}_R]}^{(i  )} - \overbrace{\mu \circ [1_R \ot \mu \circ (1_R \ot f) \circ \rho^{-1}_R]}^{(ii)} \nonumber \\
                -&\overbrace{\mu \circ (f \ot 1_R) \circ \lambda^{-1}_R \circ \mu                      }^{(iii)}   + \overbrace{\mu \circ (1_R \ot f) \circ \rho^{-1}_R \circ \mu           }^{(iv)}           \\ 
                +&\overbrace{\mu \circ [\mu \circ (f \ot 1_R) \circ \lambda^{-1}_R \ot 1_R]}^{(v  )} - \overbrace{\mu \circ [\mu \circ (1_R \ot f) \circ \rho^{-1}_R \ot 1_R]       }^{(vi)}.\nonumber
\end{align} 
First we claim that $(i)$ cancels against $(vi)$. To see this consider the following diagram, where $(i)$ is the vertical composition along left hand side, and $(vi)$ that along the right hand side.   
\begin{center}
\begin{tikzpicture}
\matrix(m)[matrix of math nodes,row sep=2.6em,column sep=2.8em,text height=1.5ex,text depth=0.25ex]
{ R\ot R             &          &   R   \ot R             \\
  R\ot (I \ot R)     &          &  (R \ot I) \ot R        \\ 
  R\ot (R \ot R)     &          &  (R \ot R) \ot R        \\ 
  R\ot            R  &          &   R \ot            R    \\ 
  R                  &          &   R                     \\ 
};
\node[above= .01cm of m-1-1,font=\scriptsize]{$(i):$};
\node[above= .01cm of m-1-3,font=\scriptsize]{$(vi):$};

\draw[ ->,font=\scriptsize](m-1-1) edge         node[left ]{$1_R\ot\lambda^{-1}_R         $} (m-2-1);
\draw[ ->,font=\scriptsize](m-2-1) edge         node[left ]{$1_R\ot(f\ot 1_R)         $} (m-3-1);
\draw[ ->,font=\scriptsize](m-3-1) edge         node[left ]{$1_R \ot \mu                  $} (m-4-1);
\draw[ ->,font=\scriptsize](m-4-1) edge         node[left ]{$\mu                              $} (m-5-1);

\draw[ ->,font=\scriptsize](m-1-3) edge         node[right]{$\rho^{-1}_R\ot 1_R           $} (m-2-3);
\draw[ ->,font=\scriptsize](m-2-3) edge         node[right]{$(1_R\ot f)\ot 1_R        $} (m-3-3);
\draw[ ->,font=\scriptsize](m-3-3) edge         node[right]{$\mu \ot 1_R                  $} (m-4-3);
\draw[ ->,font=\scriptsize](m-4-3) edge         node[right]{$\mu                              $} (m-5-3);

\draw[<- ,font=\scriptsize](m-1-1) edge[dashed] node[above]{$1_{R\ot R}                   $} (m-1-3);
\draw[<- ,font=\scriptsize](m-2-1) edge[dashed] node[above]{$\alpha_{R,1,R}                   $} (m-2-3);
\draw[<- ,font=\scriptsize](m-3-1) edge[dashed] node[above]{$\alpha_{R,R,R}                   $} (m-3-3);
\draw[<- ,font=\scriptsize](m-5-1) edge[dashed] node[above]{$1_{R}                            $} (m-5-3);

\end{tikzpicture}
\end{center}
The top square commutes from the (inverse of the) triangle identity. The middle square commutes by the naturality of $\alpha$. The bottom square commutes by the associativity relation in the multiplication rule. Hence the diagram commutes and ``$(i)+(vi)=0$'', namely $\mu \circ [1_R \ot \mu \circ (f \ot 1_R) \circ \lambda^{-1}_R]-\mu \circ [\mu \circ (1_R \ot f) \circ \rho^{-1}_R \ot 1_R] = 0$. 

Next we show that $(ii)$ cancels against $(iv)$. Consider the following diagram where $(ii)$ is the vertical composition along the left hand side and $(iv)$ that along the right hand side 
\begin{center}
\begin{tikzpicture}
\matrix(m)[matrix of math nodes,row sep=2.6em,column sep=2.8em,text height=1.5ex,text depth=0.25ex]
{ R\ot R             &                               &   R   \ot R             \\
                         &                               &   R                         \\
  R\ot (R \ot I) & (R\ot R)\ot I         &   R \ot I               \\ 
  R\ot (R \ot R) & (R\ot R)\ot R         &   R \ot R               \\ 
  R\ot            R  &                               &                             \\ 
  R                      &                               &   R                         \\ 
};
\node[above= .01cm of m-1-1,font=\scriptsize]{$(ii):$};
\node[above= .01cm of m-1-3,font=\scriptsize]{$(iv):$};

\draw[ ->,font=\scriptsize](m-1-1) edge         node[left ]{$1_R\ot\rho^{-1}_R            $} (m-3-1);
\draw[ ->,font=\scriptsize](m-3-1) edge         node[left ]{$1_R\ot(1_R\ot f)         $} (m-4-1);
\draw[ ->,font=\scriptsize](m-4-1) edge         node[left ]{$1_R \ot \mu                  $} (m-5-1);
\draw[ ->,font=\scriptsize](m-5-1) edge         node[left ]{$\mu                              $} (m-6-1);

\draw[ ->,font=\scriptsize](m-1-3) edge         node[right]{$\mu                              $} (m-2-3);
\draw[ ->,font=\scriptsize](m-2-3) edge         node[right]{$\rho^{-1}_R                      $} (m-3-3);
\draw[ ->,font=\scriptsize](m-3-3) edge         node[right]{$1_R \ot f                    $} (m-4-3);
\draw[ ->,font=\scriptsize](m-4-3) edge         node[right]{$\mu                              $} (m-6-3);

\draw[ ->,font=\scriptsize](m-1-1) edge[dashed] node[above]{$1_{R\ot R}                   $} (m-1-3);
\draw[ ->,font=\scriptsize](m-3-1) edge[dashed] node[above]{$\alpha^{-1}_{R,R,I}              $} (m-3-2);
\draw[ ->,font=\scriptsize](m-3-2) edge[dashed] node[above]{$\mu \ot 1_I                  $} (m-3-3);
\draw[ ->,font=\scriptsize](m-4-1) edge[dashed] node[above]{$\alpha^{-1}_{R,R,R}              $} (m-4-2);
\draw[ ->,font=\scriptsize](m-4-2) edge[dashed] node[above]{$\mu \ot 1_R                  $} (m-4-3);
\draw[ ->,font=\scriptsize](m-6-1) edge[dashed] node[above]{$1_{R}                            $} (m-6-3);

\draw[ ->,font=\scriptsize](m-1-3) edge[dashed] node[above]{$\rho^{-1}_{R \ot R}              $} (m-3-2);
\draw[ ->,font=\scriptsize](m-3-2) edge[dashed] node[left ]{$1_{R \ot R} \ot f                $} (m-4-2);
\end{tikzpicture}
\end{center}
The top left part of the diagram commutes by coherence and the top right part commutes since $\rho$ is natural. The middle left part commutes since $\alpha$ is natural and the middle right part of the diagram commutes by functoriality and identities. The bottom part commutes by the associative relation. Hence the diagram commutes and $(ii)+(iv)=0$, or $-\mu \circ [1_R \ot \mu \circ (1_R \ot f) \circ \rho^{-1}_R]+\mu \circ (1_R \ot f) \circ \rho^{-1}_R \circ \mu = 0$. 

Finally, we show that $(iii)+(v)=0$. In order to show this, we consider the ``$(iii)$-$(v)$''-diagram 
\begin{center}
\begin{tikzpicture}
\matrix(m)[matrix of math nodes,row sep=2.6em,column sep=2.8em,text height=1.5ex,text depth=0.25ex]
{   
    R   \ot R            &                               &  R\ot R               \\
    R                        &                               &                           \\
    I \ot R              &  I\ot(R \ot R)        & (I\ot  R)\ot R    \\ 
    R \ot R              &  R\ot(R \ot R)        & (R\ot  R)\ot R    \\ 
                             &                               &  R\ot            R    \\ 
    R                        &                               &  R                        \\ 
};
\node[above= .01cm of m-1-1,font=\scriptsize]{$(iii):$};
\node[above= .01cm of m-1-3,font=\scriptsize]{$(v):$};

\draw[ ->,font=\scriptsize](m-1-3) edge         node[right]{$\lambda^{-1}_R\ot1_R         $} (m-3-3);
\draw[ ->,font=\scriptsize](m-3-3) edge         node[right]{$(f\ot 1_R)\ot1_R         $} (m-4-3);
\draw[ ->,font=\scriptsize](m-4-3) edge         node[right]{$\mu \ot 1_R                  $} (m-5-3);
\draw[ ->,font=\scriptsize](m-5-3) edge         node[right]{$\mu                              $} (m-6-3);

\draw[ ->,font=\scriptsize](m-1-1) edge         node[left ]{$\mu                              $} (m-2-1);
\draw[ ->,font=\scriptsize](m-2-1) edge         node[left ]{$\lambda^{-1}_R                   $} (m-3-1);
\draw[ ->,font=\scriptsize](m-3-1) edge         node[left ]{$f\ot 1_R                     $} (m-4-1);
\draw[ ->,font=\scriptsize](m-4-1) edge         node[left ]{$\mu                              $} (m-6-1);

\draw[<- ,font=\scriptsize](m-1-1) edge         node[above]{$1_{R\ot R}                   $} (m-1-3);
\draw[<- ,font=\scriptsize](m-3-1) edge         node[above]{$1_I \ot \mu                  $} (m-3-2);
\draw[<- ,font=\scriptsize](m-3-2) edge         node[above]{$\alpha_{I,R,R}                   $} (m-3-3);
\draw[<- ,font=\scriptsize](m-4-1) edge         node[above]{$1_R \ot \mu                  $} (m-4-2);
\draw[<- ,font=\scriptsize](m-4-2) edge         node[above]{$\alpha_{R,R,R}                   $} (m-4-3);
\draw[<- ,font=\scriptsize](m-6-1) edge         node[above]{$1_{R}                            $} (m-6-3);

\draw[ ->,font=\scriptsize](m-1-1) edge[dashed] node[right]{$\lambda^{-1}_{R\ot R}        $} (m-3-2);
\draw[ ->,font=\scriptsize](m-3-2) edge[dashed] node[right]{$f\ot1_{R\ot R}           $} (m-4-2);
\end{tikzpicture}
\end{center}
First we observe that this diagram is the ``reflection'' of the ``$(ii)$-$(iv)$''-diagram. Again starting at the very top, the square left of the dashed $\lambda^{-1}_{R \ot R}$ commutes since $\lambda$ is natural. The square right of this dashed arrow commutes by coherence. In the middle part of the diagram, the right square commutes since $\alpha$ is natural, while the left square commutes by straightforward compositions ($-\ot-$ is a functor and properties of compositions with identities). The bottom part commutes by the associative relation for the multiplication rule. Hence $(iii)+(v)=0$, or $-\mu \circ (f \ot 1_R) \circ \lambda^{-1}_R \circ \mu + \mu \circ [ \mu \circ (f \ot 1_R) \circ \lambda^{-1}_R \ot 1_R] = 0$. We conclude that $d^1 \circ d^0=0$ from the diagrams above. 
\\\\
\textbf{The case $k\geq 1$.} We prove that $d^{k+1} \circ d^{k}=0$ for $k\geq 1$. Let $f\in\Hom_{\fC}(R^{\ot k},R)$ and recall that 
\begin{align*}
d^{k+1}(?)={ {\mu} \circ {(1_R \ot ?)} \circ \alpha^{0,2}_k }  
      +\sum^{k+1}_{i=1}(-1)^i[{?} \circ {(\alpha^{i-1,1}_{k+1})^{-1} \circ {\mu^{i-1}_{k+2}} \circ \alpha^{i-1,2}_{k+2}}] 
      +(-1)^{k+2}[{\mu} \circ {(? \ot 1_R)}]. 
\end{align*}
Since 
\begin{align*}
d^k(f)={ {\mu} \circ {(1_R \ot f)} \circ {\alpha^{0,1}_{k+1}}}  
      +\sum^{k}_{i=1}(-1)^i[{f} \circ (\alpha^{i-1,1}_k)^{-1} \circ {\mu^{i-1}_{k+1}} \circ \alpha^{i-1,2}_{k+1}] 
      +(-1)^{k+1}[{\mu} \circ {(f \ot 1_R)}]\nonumber \\ 
\end{align*}
we then obtain 
\begin{align}
(d^{k+1} \circ d^k)(f)=\nonumber\\ 
             &   \overbrace{\mu \circ [ 1_R \ot \mu \circ (1_R \ot f) \circ \alpha^{0,1}_{k+1}] \circ \alpha^{0,1}_{k+2}                                                                                                             }^{(i    )}       \nonumber\\ 
             &+  \overbrace{\mu \circ \left[1_R\ot \left(\sum^{k}_{i=1}(-1)^i{f} \circ {(\alpha^{i-1,1}_k)^{-1}} \circ {\mu^{i-1}_{k+1}} \circ {\alpha^{i-1,2}_{k+1}}\right)\right] \circ \alpha^{0,1}_{k+2}                                       }^{(ii   )}\nonumber\\ 
             &+  \overbrace{(-1)^{k+1}\mu \circ [ 1_R \ot \mu \circ (f \ot 1_R) ]\circ \alpha^{0,1}_{k+2}                                                                                                                    }^{(iii  )} \nonumber\\
             &+  \overbrace{\sum^{k+1}_{i=1}(-1)^i[\mu \circ (1_R \ot f) \circ \alpha^{0,1}_{k+1}](\alpha^{i-1,1}_{k+1})^{-1} \circ \mu^{i-1}_{k+2} \circ \alpha^{i-1,2}_{k+2}                                                            }^{(iv   )}\nonumber\\ 
             &+  \overbrace{\sum^{k+1}_{i=1}(-1)^{i}\left[\sum^{k}_{j=1}(-1)^{j} f \circ (\alpha^{j-1,1}_{k})^{-1} \circ \mu^{j-1}_{k+1} \circ \alpha^{j-1,2}_{k+1}\right](\alpha^{i-1,1}_{k+1})^{-1} \circ \mu^{i-1}_{k+2} \circ \alpha^{i-1,2}_{k+2} }^{(v    )}\nonumber\\
             &+  \overbrace{\sum^{k+1}_{i=1}(-1)^{i}(-1)^{k+1}\mu \circ (f\ot 1_R) \circ (\alpha^{i-1,1}_{k+1})^{-1} \circ \mu^{i-1}_{k+2} \circ \alpha^{i-1,2}_{k+2}                                                                          }^{(vi   )}\nonumber\\ 
             &+  \overbrace{(-1)^{k+2}\mu \circ ([ \mu \circ (1_R \ot f) \circ \alpha^{0,1}_{k+1}] \ot 1_R )                                                                                                                    }^{(vii  )}\nonumber\\
             &+  \overbrace{(-1)^{k+2}\mu \circ \bigg[\sum^k_{i=1}(-1)^i {f} \circ {(\alpha^{i-1,1}_k)^{-1}} \circ {\mu^{i-1}_{k+1}} \circ {\alpha^{i-1,2}_{k+1}} \ot 1_R \bigg]                                                        }^{(viii )}\nonumber\\
             &+  \overbrace{(-1)^{k+2}(-1)^{k+1} \mu \circ [ \mu \circ (f \ot 1_R) \ot 1_R ]                                                                                                                           }^{(ix   )}.        
\end{align}

Before we prove that this expression vanishes, we introduce the following notation. Let $(??)=\sum_{i=1}^k\xi_i$ be one of the nine sums above. For $0 \leq n \leq k$, by $(??)_{n}$ we denote the $n$th term in expression, that is $(??)_n=\xi_n$. Furthermore $(??)^{n}$ denotes the sum $(??)$ but without the $n$th term, that is $(??)^{n}=\sum_{i=1}^{n-1}\xi_i+\sum_{i=n+1}^{k}\xi_i$. 

First we associate $(iii)$ with $(vii)$. We use that $\ot$ is a bifunctor and hence commutes with compositions to rewrite $(iii)$ and $(vii)$ slightly  
\begin{align*}
(iii)=&(-1)^{k+1}\mu \circ [ 1_R \ot \mu \circ (f \ot 1_R) ]\circ \alpha^{0,1}_{k+2}\\
     =&(-1)^{k+1}\mu \circ (1_R\ot\mu)\circ (1_R\ot[f\ot1_R])\circ \alpha^{0,1}_{k+2}\\ 
(vii)=&(-1)^{k+2}\mu \circ ([ \mu \circ (1_R\ot f) \circ \alpha^{0,1}_{k+1}] \ot 1_R )\\
     =&(-1)^{k+2}\mu \circ (\mu\ot1_R) \circ ([1_R\ot f]\ot1_R) \circ (\alpha^{0,1}_{k+1}\ot1_R). 
\end{align*} 
Now we use these expressions to construct the following diagram: 
\begin{center}
\begin{tikzpicture}
\matrix(m)[matrix of math nodes,row sep=2.6em,column sep=2.8em,text height=1.5ex,text depth=0.25ex]
{   
    R^{\ot(k+2)}      &       &    R^{\ot(k+2)}          \\
    R\ot R^{\ot(k+1)} &       &    (R\ot R^{\ot k})\ot R \\
    R\ot(R\ot R)      &       &    (R\ot R) \ot R        \\
    R\ot R            &       &    R\ot R                \\
    R                 &       &    R                     \\
};
\node[above= .01cm of m-1-1,font=\scriptsize]{$(iii):$};
\node[above= .01cm of m-1-3,font=\scriptsize]{$(vii):$};

\draw[ ->,font=\scriptsize](m-1-1) edge         node[left ]{$ \alpha^{0,1}_{k+2}          $} (m-2-1);
\draw[ ->,font=\scriptsize](m-2-1) edge         node[left ]{$ 1_R\ot[f\ot1_R]             $} (m-3-1);
\draw[ ->,font=\scriptsize](m-3-1) edge         node[left ]{$ 1_R \ot \mu                 $} (m-4-1);
\draw[ ->,font=\scriptsize](m-4-1) edge         node[left ]{$ \mu                         $} (m-5-1);

\draw[ ->,font=\scriptsize](m-1-3) edge         node[right]{$ \alpha^{0,1}_{k+1}\ot1_R    $} (m-2-3);
\draw[ ->,font=\scriptsize](m-2-3) edge         node[right]{$ (1_R\ot f)\ot 1_R           $} (m-3-3);
\draw[ ->,font=\scriptsize](m-3-3) edge         node[right]{$ \mu \ot 1_R                 $} (m-4-3);
\draw[ ->,font=\scriptsize](m-4-3) edge         node[right]{$ \mu                         $} (m-5-3);

\draw[<- ,font=\scriptsize](m-1-1) edge[dashed] node[above]{$ 1_{R^{\ot (k+1)}}           $} (m-1-3);
\draw[<- ,font=\scriptsize](m-2-1) edge[dashed] node[above]{$ \alpha_{R,R^{\ot k},R}      $} (m-2-3);
\draw[<- ,font=\scriptsize](m-3-1) edge[dashed] node[above]{$ \alpha_{R,R,R}              $} (m-3-3);
\draw[<- ,font=\scriptsize](m-5-1) edge[dashed] node[above]{$ 1_R                         $} (m-5-3);
\end{tikzpicture}
\end{center}
The top square commutes by the coherence theorem. The middle part commutes since $\alpha$ is a natural transformation. The bottom part commutes by the associativity relation of the multiplication rule. Hence the diagram as a whole commutes and we conclude that $(iii)+(vii)=0$. 

Next we associate $(i)$ to $(iv)_1$. Again we rewrite sightly
\begin{align*}
(i)   &=\mu \circ [ 1_R \ot \mu (1_R \ot f) \circ \alpha^{0,1}_{k+1} ]\circ \alpha^{0,1}_{k+2}\\
      &=\mu \circ (1_R\ot\mu) \circ (1_R\ot[1_R\ot f]) \circ (1_R\ot\alpha^{0,1}_{k+1}) \circ \alpha^{0,1}_{k+2} \\
(iv)_1&=(-1)[\mu \circ (1_R \ot f) \circ \alpha^{0,1}_{k+1} ] \circ (\alpha^{0,1}_{k+1})^{-1} \circ \mu^{0}_{k+2} \circ \alpha^{0,2}_{k+2}\\
      &=(-1) \mu \circ (1_R \ot f) \circ \mu^{0}_{k+2} \circ \alpha^{0,2}_{k+2}
\end{align*}
Then we organise these expressions in the following diagram: 
\begin{center}
\begin{tikzpicture}
\matrix(m)[matrix of math nodes,row sep=2.6em,column sep=2.8em,text height=1.5ex,text depth=0.25ex]
{   
    R^{\ot(k+2)}                            &                             & R^{\ot (k+2)}                               \\
    R\ot R^{\ot (k+1)}                  &                             &                                             \\
    R\ot( R\ot R^{\ot k})           &                             & ( R\ot  R)\ot R^{\ot k}             \\
    R\ot (R \ot R)                  &                             &   R\ot  R^{\ot k}                       \\
    R\ot R                              &  (R \ot R) \ot R    &   R\ot  R                               \\ 
    R                                       &                             &   R                                         \\ 
};
\node[above= .01cm of m-1-1,font=\scriptsize]{$(i ):$};
\node[above= .01cm of m-1-3,font=\scriptsize]{$(iv):$};

\draw[ ->,font=\scriptsize](m-1-1) edge         node[left ]{$ \alpha^{0,1}_{k+2}              $} (m-2-1);
\draw[ ->,font=\scriptsize](m-2-1) edge         node[left ]{$ 1_R\ot\alpha^{0,1}_{k+1}    $} (m-3-1);
\draw[ ->,font=\scriptsize](m-3-1) edge         node[left ]{$ 1_R\ot(1_R\ot f)        $} (m-4-1);
\draw[ ->,font=\scriptsize](m-4-1) edge         node[left ]{$ 1_R\ot\mu                   $} (m-5-1);
\draw[ ->,font=\scriptsize](m-5-1) edge         node[left ]{$ \mu                             $} (m-6-1);

\draw[ ->,font=\scriptsize](m-1-3) edge         node[right]{$ \alpha^{0,1}_{k+2}              $} (m-3-3);
\draw[ ->,font=\scriptsize](m-3-3) edge         node[right]{$ \mu\ot 1_{R^{\ot k}}    $} (m-4-3);
\draw[ ->,font=\scriptsize](m-4-3) edge         node[right]{$ 1_R \ot f                   $} (m-5-3);
\draw[ ->,font=\scriptsize](m-5-3) edge         node[right]{$ \mu                             $} (m-6-3);

\draw[ ->,font=\scriptsize](m-1-1) edge[dashed] node[above]{$1_{R^{\ot (k+2)}}                $} (m-1-3);
\draw[ ->,font=\scriptsize](m-3-1) edge[dashed] node[above]{$\alpha^{-1}_{R,R,R^{\ot k}}  $} (m-3-3);
\draw[ ->,font=\scriptsize](m-4-1) edge[dashed] node[above]{$\alpha^{-1}_{R,R,R}              $} (m-5-2);
\draw[ ->,font=\scriptsize](m-3-3) edge[dashed] node[left ]{$1_{R\ot R}\ot f          $} (m-5-2);
\draw[ ->,font=\scriptsize](m-5-2) edge[dashed] node[above]{$\mu\ot1_R                    $} (m-5-3);
\draw[ ->,font=\scriptsize](m-6-1) edge[dashed] node[above]{$1_R                              $} (m-6-3);
\end{tikzpicture}
\end{center}
The top square commutes by the coherence theorem. For the middle part, the left square commutes by the naturality of $\alpha$, and the right middle square commutes since $\ot$ is a functor. Finally, the bottom part commutes by the associativity relation. This proves that the diagram commutes and $(i)+(iv)_1=0$. 

The next objective is to prove that $(vi)_{k+1}+(ix)=0$. We use that $(\alpha^{k,1}_{k+1})^{-1}=1_{R^{\ot (k+1)}}$ and $(-1)^{k+1}(-1)^{k+1}=1$ to rewrite 
\begin{align*}
(vi)_{k+1}&=(-1)^{k+1}(-1)^{k+1}\mu \circ (f\ot1_R) \circ (\alpha^{k,1}_{k+1})^{-1} \circ \mu^{k}_{k+2} \circ \alpha^{k,2}_{k+2} \\ 
          &=\mu \circ (f\ot1_R) \circ (1_{R^{\ot k}}\ot\mu) \circ \alpha^{k,2}_{k+2}\\
(ix)      &=(-1)^{k+2}(-1)^{k+1}\mu \circ [\mu \circ (f\ot1_R)\ot1_R]\\
          &=-\mu \circ (\mu\ot1_R) \circ ([f\ot1_R]\ot1_R) 
\end{align*}
Now consider the diagram 
\begin{center}
\begin{tikzpicture}
\matrix(m)[matrix of math nodes,row sep=2.6em,column sep=2.8em,text height=1.5ex,text depth=0.25ex]
{   
    R^{\ot(k+2)}                                                        &                             & R^{\ot (k+2)}                                                          \\
    R^{\ot k}\ot(R\ot R)                                                &                             &                                                                        \\
    R^{\ot k }\ot R                                                     &                             & ( R\ot  R)\ot R                                                \\
    R^{\ot (k+1)}                                                       & R \ot (R\ot R )             &                                                                        \\
    R\ot R                                                          &                             &   R\ot  R                                                          \\ 
    R                                                                   &                             &   R                                                                    \\ 
};
\node[above= .01cm of m-1-1,font=\scriptsize]{$(vi)_{k+1}:$};
\node[above= .01cm of m-1-3,font=\scriptsize]{$(ix):$};

\draw[ ->,font=\scriptsize](m-1-1) edge         node[left ]{$ \alpha^{k,2}_{k+2}              $} (m-2-1);
\draw[ ->,font=\scriptsize](m-2-1) edge         node[left ]{$ 1_{R^{\ot k}}\ot \mu            $} (m-3-1);
\draw[ ->,font=\scriptsize](m-3-1) edge         node[left ]{$ 1_{R^{\ot (k+1)}}               $} (m-4-1);
\draw[ ->,font=\scriptsize](m-4-1) edge         node[left ]{$ f\ot 1_R                        $} (m-5-1);
\draw[ ->,font=\scriptsize](m-5-1) edge         node[left ]{$ \mu                             $} (m-6-1);
\draw[ ->,font=\scriptsize](m-1-3) edge         node[right]{$ [f\ot 1_R]\ot 1_R               $} (m-3-3);
\draw[ ->,font=\scriptsize](m-3-3) edge         node[right]{$ \mu\ot 1_{R}                $} (m-5-3);
\draw[ ->,font=\scriptsize](m-5-3) edge         node[right]{$ \mu                             $} (m-6-3);
\draw[ ->,font=\scriptsize](m-1-3) edge[dashed] node[above]{$1_{R^{\ot (k+2)}}                $} (m-1-1);
\draw[ ->,font=\scriptsize](m-1-3) edge[dashed] node[below,sloped]{$\alpha_{R^{\ot k},R,R}           $} (m-2-1);
\draw[ ->,font=\scriptsize](m-2-1) edge[dashed] node[above,sloped]{$f\ot1_{R^{\ot2}}                 $} (m-4-2);
\draw[ ->,font=\scriptsize](m-3-3) edge[dashed] node[above,sloped]{$\alpha_{R,R,R}                   $} (m-4-2);
\draw[ ->,font=\scriptsize](m-4-2) edge[dashed] node[above,sloped]{$1_R\ot\mu                        $} (m-5-1);
\draw[ ->,font=\scriptsize](m-6-3) edge[dashed] node[above,sloped]{$1_R                              $} (m-6-1);
\end{tikzpicture}
\end{center}
The top triangle commutes by the coherence theorem. The middle part of the diagram consists of two squares, 
\begin{center}
\begin{tikzpicture}
\matrix(m)[matrix of math nodes,row sep=2.6em,column sep=2.8em,text height=1.5ex,text depth=0.25ex]
{   
    R^{\ot k }\ot R           & R^{\ot k}\ot(R\ot R)        & R^{\ot (k+2)}                   \\
    R\ot R                & R \ot (R\ot R )             & ( R\ot  R)\ot R         \\ 
};
\draw[ ->,font=\scriptsize](m-1-3) edge[dashed] node[above]{$ \alpha_{R^{\ot k},R,R}          $} (m-1-2);
\draw[ ->,font=\scriptsize](m-1-2) edge         node[above]{$ 1_{R^{\ot k}}\ot \mu            $} (m-1-1);

\draw[ ->,font=\scriptsize](m-1-1) edge         node[left ]{$ f\ot 1_R                        $} (m-2-1);
\draw[ ->,font=\scriptsize](m-1-3) edge         node[left ]{$ [f\ot 1_R]\ot 1_R               $} (m-2-3);
\draw[ ->,font=\scriptsize](m-1-2) edge[dashed] node[left ]{$f\ot1_{R^{\ot2}}                 $} (m-2-2);

\draw[ ->,font=\scriptsize](m-2-3) edge[dashed] node[above]{$\alpha_{R,R,R}                   $} (m-2-2);
\draw[ ->,font=\scriptsize](m-2-2) edge[dashed] node[above]{$1_R\ot\mu                        $} (m-2-1);
\end{tikzpicture}
\end{center}
where the right square commutes by the naturality of $\alpha$. The left square clearly commutes, hence the middle part of the diagram commutes. The bottom part commutes by the multiplicative associativity relation, and we conclude that $(vi)_{n+1}+(ix)=0$. 

Next we prove that the sums $(ii)+(iv)^1=0$. We do this by checking that $(ii)_{i}+(iv)_{i+1}^1=0$ for $0\leq i \leq k$. Recall that  
\begin{align*}
(ii)_{i}    &=\mu \circ [1_R\ot ((-1)^i{f} \circ {(\alpha^{i-1,1}_k)^{-1}} \circ {\mu^{i-1}_{k+1}} \circ {\alpha^{i-1,2}_{k+1}})] \circ \alpha^{0,1}_{k+2}\\
(iv)_{i+1}^1&=(-1)^{i+1}[\mu \circ (1_R \ot f) \circ \alpha^{0,1}_{k+1} ] \circ (\alpha^{i,1}_{k+1})^{-1} \circ \mu^{i}_{k+2} \circ \alpha^{i,2}_{k+2} 
\end{align*}
and consider that diagram 
\begin{center}
\begin{tikzpicture}
\matrix(m)[matrix of math nodes,row sep=2.6em,column sep=2.8em,text height=1.5ex,text depth=0.25ex]
{   
 R^{\ot(k+2)}                                      & R^{\ot (k+2)}                                   \\
 R\ot R^{\ot k+1}                                  &                                                 \\
 R\ot([R^{\ot(i-1)}\ot R^{\ot 2}]\ot R^{\ot(k-i)}) & ( R^{\ot i}\ot  R^{\ot 2})\ot R^{\ot (k-i)} \\
 R\ot([R^{\ot(i-1)}\ot R        ]\ot R^{\ot(k-i)}) & ( R^{\ot i}\ot  R        )\ot R^{\ot (k-i)} \\
                                                   &   R^{\ot(k+1)}                                  \\
 R\ot R^{\ot k}                                    &   R\ot R^{\ot k}                                \\
 R\ot R                                        &   R\ot  R                                   \\ 
 R                                                 &   R                                             \\ 
};
\node[above= .01cm of m-1-1,font=\scriptsize]{$(ii)_{i}:    $};
\node[above= .01cm of m-1-2,font=\scriptsize]{$(iv)_{i+1}^1:$};
\draw[ ->,font=\scriptsize](m-1-1) edge         node[left ]{$\alpha^{0,1}_{k+2}                                  $} (m-2-1);
\draw[ ->,font=\scriptsize](m-2-1) edge         node[left ]{$1_R\ot\alpha^{i-1,2}_{k+1}                          $} (m-3-1);
\draw[ ->,font=\scriptsize](m-3-1) edge         node[left ]{$1_R\ot\mu^{i-1}_{k+1}                               $} (m-4-1);
\draw[ ->,font=\scriptsize](m-4-1) edge         node[left ]{$1_R\ot(\alpha^{i-1,1}_k)^{-1}                       $} (m-6-1);
\draw[ ->,font=\scriptsize](m-6-1) edge         node[left ]{$1_R\ot f                                            $} (m-7-1);
\draw[ ->,font=\scriptsize](m-7-1) edge         node[left ]{$\mu                                                 $} (m-8-1);
\draw[ ->,font=\scriptsize](m-1-2) edge         node[left ]{$\alpha^{i,2}_{k+2}                                  $} (m-3-2);
\draw[ ->,font=\scriptsize](m-3-2) edge         node[left ]{$\mu^{i}_{k+2}                                       $} (m-4-2);
\draw[ ->,font=\scriptsize](m-4-2) edge         node[left ]{$(\alpha^{i,1}_{k+1})^{-1}                           $} (m-5-2);
\draw[ ->,font=\scriptsize](m-5-2) edge         node[left ]{$\alpha^{0,1}_{k+1}                                  $} (m-6-2);
\draw[ ->,font=\scriptsize](m-6-2) edge         node[left ]{$1_R\ot f                                            $} (m-7-2);
\draw[ ->,font=\scriptsize](m-7-2) edge         node[left ]{$\mu                                                 $} (m-8-2);
\draw[ ->,font=\scriptsize](m-1-1) edge[dashed] node[above]{$1_{R^{\ot (k+2)}}                                   $} (m-1-2);
\draw[ ->,font=\scriptsize](m-3-1) edge[dashed] node[above]{                                                      } (m-3-2);
\draw[ ->,font=\scriptsize](m-4-1) edge[dashed] node[above]{                                                      } (m-4-2);
\draw[ ->,font=\scriptsize](m-6-1) edge[dashed] node[above]{$1_{R\ot R^{\ot k}}                                  $} (m-6-2);
\draw[ ->,font=\scriptsize](m-7-1) edge[dashed] node[above]{$1_{R\ot R}                                          $} (m-7-2);
\draw[ ->,font=\scriptsize](m-8-1) edge[dashed] node[above]{$1_{R}                                               $} (m-8-2);
\end{tikzpicture}
\end{center}
The unlabeled horizontal arrows are compositions of associators appropriate to the setting. The top part commutes by the coherence theorem and the second square commutes by the naturality of the associator. The middle square commutes by the coherence theorem again, while the bottom two squares clearly commute. The diagram shows that $(ii)_{i}+(iv)_{i+1}^1=0$ for all $0\leq i \leq k$, hence $(ii)+(iv)^1$ for the full sums. 

Next we consider the ``reflected'' version of the previous identification. We want to show that $(vi)^{k+1}+(viii)=0$, that is 
\begin{align*}
0=(vi)^{k+1}+(viii)
=\sum^{k}_{i=1}(-1)^{i}(-1)^{k+1} \mu \circ (f\ot 1_R) \circ (\alpha^{i-1,1}_{k+1})^{-1} \circ \mu^{i-1}_{k+2} \circ \alpha^{i-1,2}_{k+2}+\\
 (-1)^{k+2} \mu \circ \left[\sum^k_{i=1}(-1)^i {f} \circ {(\alpha^{i-1,1}_k)^{-1}} \circ {\mu^{i-1}_{k+1}} \circ {\alpha^{i-1,2}_{k+1}} \ot 1_R \right]\\
=\sum^{k}_{i=1}\left( \mu \circ (f\ot 1_R) \circ (\alpha^{i-1,1}_{k+1})^{-1} \circ \mu^{i-1}_{k+2} \circ \alpha^{i-1,2}_{k+2}-
 \mu \circ \left[{f} \circ {(\alpha^{i-1,1}_k)^{-1}} \circ {\mu^{i-1}_{k+1}} \circ {\alpha^{i-1,2}_{k+1}} \ot 1_R \right]\right)
\end{align*}
In this last sum we check that each term vanishes, i.e.\ 
\[\mu \circ (f\ot 1_R) \circ (\alpha^{i-1,1}_{k+1})^{-1} \circ \mu^{i-1}_{k+2} \circ \alpha^{i-1,2}_{k+2}-\mu \circ \left[{f} {(\alpha^{i-1,1}_k)^{-1}} \circ {\mu^{i-1}_{k+1}} \circ {\alpha^{i-1,2}_{k+1}} \ot 1_R \right]=0, \] 
by considering the following diagram 
\begin{center}
\begin{tikzpicture}
\matrix(m)[matrix of math nodes,row sep=2.6em,column sep=2.8em,text height=1.5ex,text depth=0.25ex]
{   
    R^{\ot(k+2)}                                           &    &    &  R^{\ot (k+2)}                                                         \\
   (R^{\ot(i-1)}\ot R^{\ot 2})\ot R^{\ot(k-i+1)}           &    &    &((R^{\ot(i-1)}\ot R^{\ot 2})\ot R^{\ot(k-i)})\ot R                      \\
   (R^{\ot(i-1)}\ot R        )\ot R^{\ot(k-i+1)}           &    &    &((R^{\ot(i-1)}\ot R        )\ot R^{\ot(k-i)})\ot R                      \\
    R^{\ot(k+1)}                                           &    &    &  R^{\ot(k+1)}                                                          \\
    R^{\ot2}                                               &    &    &  R^{\ot2}                                                              \\
    R                                                      &    &    &  R                                                                     \\
};
\node[above= .01cm of m-1-1,font=\scriptsize]{$(vi)_{i}^{k+1}:$};
\node[above= .01cm of m-1-4,font=\scriptsize]{$(viii)_{i}:$};
\draw[ ->,font=\scriptsize](m-1-1) edge         node[left ]{$\alpha^{i-1,2}_{k+2}                                       $} (m-2-1);
\draw[ ->,font=\scriptsize](m-2-1) edge         node[left ]{$\mu^{i-1}_{k+2}                                            $} (m-3-1);
\draw[ ->,font=\scriptsize](m-3-1) edge         node[left ]{$(\alpha^{i-1,1}_{k+1})^{-1}                                $} (m-4-1); 
\draw[ ->,font=\scriptsize](m-4-1) edge         node[left ]{$f\ot1_R                                                    $} (m-5-1);
\draw[ ->,font=\scriptsize](m-5-1) edge         node[left ]{$\mu                                                        $} (m-6-1);
\draw[ ->,font=\scriptsize](m-1-4) edge         node[right]{$\alpha^{i-1,2}_{k+1}           \ot1_{R}                    $} (m-2-4);
\draw[ ->,font=\scriptsize](m-2-4) edge         node[right]{$\mu^{i-1}_{k+1}                \ot1_{R}                    $} (m-3-4);
\draw[ ->,font=\scriptsize](m-3-4) edge         node[right]{$(\alpha^{i-1,1}_{k})^{-1}      \ot1_{R}                    $} (m-4-4); 
\draw[ ->,font=\scriptsize](m-4-4) edge         node[right]{$f\ot1_R                                                    $} (m-5-4);
\draw[ ->,font=\scriptsize](m-5-4) edge         node[right]{$\mu                                                        $} (m-6-4);
\draw[<- ,font=\scriptsize](m-1-1) edge[dashed] node[above]{$1_{R^{\ot(k+2)}}                                           $} (m-1-4);
\draw[<- ,font=\scriptsize](m-2-1) edge[dashed] node[above]{$                                                           $} (m-2-4);
\draw[<- ,font=\scriptsize](m-3-1) edge[dashed] node[above]{$                                                           $} (m-3-4);
\draw[<- ,font=\scriptsize](m-4-1) edge[dashed] node[above]{$1_{R^{\ot(k+1)}}                                           $} (m-4-4);
\draw[<- ,font=\scriptsize](m-5-1) edge[dashed] node[above]{$1_{R^{\ot2}}                                               $} (m-5-4);
\draw[<- ,font=\scriptsize](m-6-1) edge[dashed] node[above]{$1_R                                                        $} (m-6-4);
\end{tikzpicture}
\end{center}
The unlabeled dashed arrows are associativity relations. The top square commutes by the coherence theorem. The second square from the top commutes since the associator is natural. The middle square commutes again by the coherence theorem. While the two bottom squares commute simply by successive compositions with identities. Hence $(vi)^{k+1}+(viii)=0$. 

Finally,  we show that $(v)=0$, that is
\begin{align*}
(v)&=\sum^{k+1}_{i=1}(-1)^{i}\left[\sum^{k}_{j=1}(-1)^{j}f \circ (\alpha^{j-1,1}_{k})^{-1} \circ \mu^{j-1}_{k+1} \circ \alpha^{j-1,2}_{k+1}\right] \circ (\alpha^{i-1,1}_{k+1})^{-1} \circ \mu^{i-1}_{k+2} \circ \alpha^{i-1,2}_{k+2}\\
   &=\sum^{k+1}_{i=1}\left[\sum^{k}_{j=1}b_j\right]a_i=\sum^{k+1}_{i=1}\sum^{k}_{j=1}b_ja_i=0, 
\end{align*}
where $b_j=(-1)^{j}f \circ (\alpha^{j-1,1}_{k})^{-1} \circ \mu^{j-1}_{k+1} \circ \alpha^{j-1,2}_{k+1}$ and $a_i=(\alpha^{i-1,1}_{k+1})^{-1} \circ \mu^{i-1}_{k+2} \circ \alpha^{i-1,2}_{k+2}$. To do this, we show that $b_ja_i+b_ia_{j+1}=0$ for $1\leq i,j \leq k$. Clearly $b_{j}a_{i}$ and $b_{i}a_{j+1}$ have opposite signs, namely $(-1)^{i+j}$ and $(-1)^{i+j+1}$, respectively. Observe also that every term in the sum $(v)$ fits with this description, so if the claim is true, then $(v)=0$. First we check the case when $i=j$, i.e.\ we show that $b_{i}a_{i}+b_{i}a_{i+1}=0$. We construct the following diagram, where $b_ia_i$ is the composition along the outer left hand side, and $b_ia_{i+1}$ is that of the outer right hand side. 
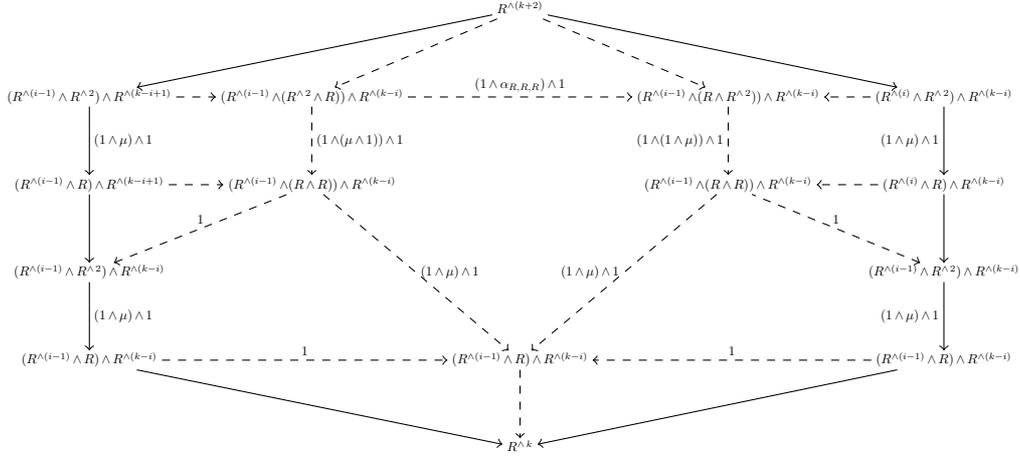
\begin{center}
\begin{tikzpicture}[scale=0.5, every node/.style={scale=0.5}]
\matrix(m)[matrix of math nodes,row sep=2.6em,column sep=1.5em,text height=1.5ex,text depth=0.25ex]
{   
                                             &                               & R^{\ot(k+2)}&                                   &                                        \\
 (R^{\ot(i-1)}\ot R^{\ot2})\ot R^{\ot(k-i+1)}&(R^{\ot(i-1)}\ot (R^{\ot2}\ot R))\ot R^{\ot(k-i)}&&(R^{\ot(i-1)}\ot (R\ot R^{\ot2}))\ot R^{\ot(k-i)}&(R^{\ot(i  )}\ot R^{\ot2})\ot R^{\ot(k-i)}\\
 (R^{\ot(i-1)}\ot R       )\ot R^{\ot(k-i+1)}&(R^{\ot(i-1)}\ot (R       \ot R))\ot R^{\ot(k-i)}&&(R^{\ot(i-1)}\ot (R\ot R       ))\ot R^{\ot(k-i)}&(R^{\ot(i  )}\ot R       )\ot R^{\ot(k-i)}\\
 (R^{\ot(i-1)}\ot R^{\ot2})\ot R^{\ot(k-i  )}&                                                 &&                                                 &(R^{\ot(i-1)}\ot R^{\ot2})\ot R^{\ot(k-i)}\\
 (R^{\ot(i-1)}\ot R       )\ot R^{\ot(k-i  )}&                               &(R^{\ot(i-1)}\ot R)\ot R^{\ot(k-i)}&                                &(R^{\ot(i-1)}\ot R       )\ot R^{\ot(k-i)}\\
                                             &                               & R^{\ot{k}}   &                                    &                                         \\
};
\draw[ ->](m-1-3) edge         node[left ]{$                          $} (m-2-1);
\draw[ ->](m-1-3) edge[dashed] node[left ]{$                          $} (m-2-2);
\draw[ ->](m-1-3) edge[dashed] node[left ]{$                          $} (m-2-4); 
\draw[ ->](m-1-3) edge         node[left ]{$                          $} (m-2-5);
\draw[ ->](m-2-1) edge         node[right]{$(1\ot\mu)\ot1             $} (m-3-1);
\draw[ ->](m-2-2) edge[dashed] node[right]{$(1\ot(\mu\ot1))\ot1       $} (m-3-2);
\draw[ ->](m-2-4) edge[dashed] node[left ]{$(1\ot(1\ot\mu))\ot1       $} (m-3-4);
\draw[ ->](m-2-5) edge         node[left ]{$(1\ot\mu)\ot1             $} (m-3-5);
\draw[ ->](m-2-1) edge[dashed] node[left ]{$                          $} (m-2-2);
\draw[ ->](m-3-1) edge[dashed] node[left ]{$                          $} (m-3-2);
\draw[ ->](m-2-2) edge[dashed] node[above]{$(1\ot\alpha_{R,R,R})\ot1  $} (m-2-4);
\draw[ ->](m-2-5) edge[dashed] node[left ]{$                          $} (m-2-4);
\draw[ ->](m-3-5) edge[dashed] node[left ]{$                          $} (m-3-4);
\draw[ ->](m-3-1) edge         node[left ]{$                          $} (m-4-1);
\draw[ ->](m-3-5) edge         node[left ]{$                          $} (m-4-5);
\draw[<- ](m-4-1) edge[dashed] node[above]{$1                         $} (m-3-2);
\draw[<- ](m-4-5) edge[dashed] node[above]{$1                         $} (m-3-4);
\draw[ ->](m-3-2) edge[dashed] node[right]{$(1\ot\mu)\ot1             $} (m-5-3);
\draw[ ->](m-3-4) edge[dashed] node[left ]{$(1\ot\mu)\ot1             $} (m-5-3);
\draw[ ->](m-4-1) edge         node[right]{$(1\ot\mu)\ot1             $} (m-5-1);
\draw[ ->](m-4-5) edge         node[left ]{$(1\ot\mu)\ot1             $} (m-5-5);
\draw[ ->](m-5-1) edge[dashed] node[above]{$1                         $} (m-5-3);
\draw[ ->](m-5-5) edge[dashed] node[above]{$1                         $} (m-5-3);
\draw[ ->](m-5-1) edge         node[above]{$                          $} (m-6-3);
\draw[ ->](m-5-3) edge[dashed] node[left ]{$                          $} (m-6-3);
\draw[ ->](m-5-5) edge         node[above]{$                          $} (m-6-3);
\end{tikzpicture}
\end{center}
The unlabeled arrows are the associativity relations. The middle part of the diagram commutes by the associativity relation for the multiplication rule. The rest of the diagram commutes by naturality, identities and the coherence theorem. For the remaining case when $i<j$ we have the following diagram, again where $b_{j}a_{i}$ is the composition along the outer left hand side and $b_{i}a_{j+1}$ is that of the outer right hand side.  
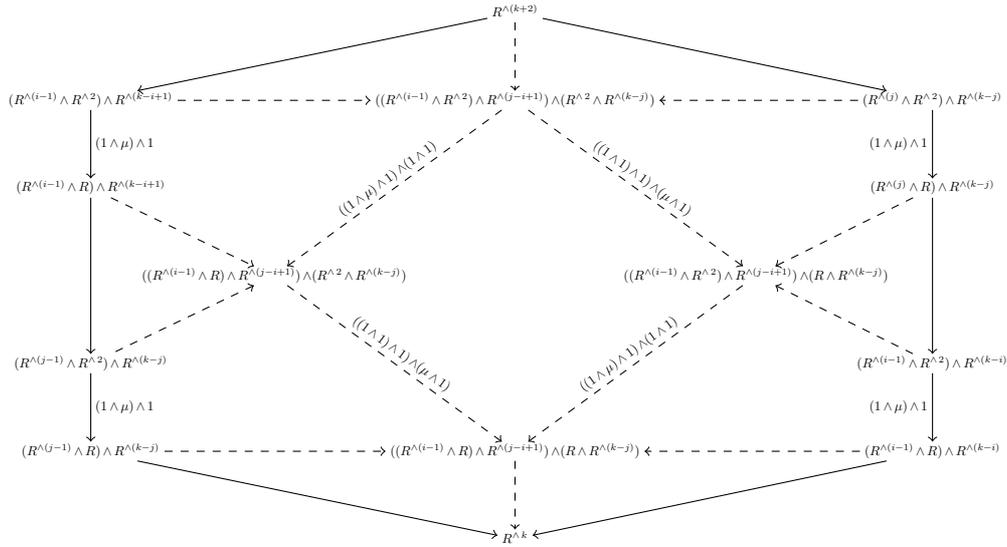
\begin{center}
\begin{tikzpicture}[scale=0.5, every node/.style={scale=0.5}]
\matrix(m)[matrix of math nodes,row sep=2.6em,column sep=-1.5em,text height=1.5ex,text depth=0.25ex]
{   
                                            &                               & R^{\ot(k+2)}&                                   &                                    \\
(R^{\ot(i-1)}\ot R^{\ot2})\ot R^{\ot(k-i+1)}&&((R^{\ot(i-1)}\ot R^{\ot2})\ot R^{\ot(j-i+1)})\ot(R^{\ot2}\ot R^{\ot(k-j)})&&(R^{\ot(j)}\ot R^{\ot2})\ot R^{\ot(k-j)}\\
(R^{\ot(i-1)}\ot R)\ot R^{\ot(k-i+1)}       &&                                                                           &&(R^{\ot(j)}\ot R       )\ot R^{\ot(k-j)}\\
         &((R^{\ot(i-1)}\ot R)\ot R^{\ot(j-i+1)})\ot(R^{\ot2}\ot R^{\ot(k-j)})&&((R^{\ot(i-1)}\ot R^{\ot2})\ot R^{\ot(j-i+1)})\ot(R\ot R^{\ot(k-j)})&  \\
(R^{\ot(j-1)}\ot R^{\ot2})\ot R^{\ot(k-j  )}&                                                 &&                                                 &(R^{\ot(i-1)}\ot R^{\ot2})\ot R^{\ot(k-i)}\\
(R^{\ot(j-1)}\ot R       )\ot R^{\ot(k-j  )}&&((R^{\ot(i-1)}\ot R       )\ot R^{\ot(j-i+1)})\ot(R\ot R^{\ot(k-j)})&&(R^{\ot(i-1)}\ot R       )\ot R^{\ot(k-i)}\\
                                            &                               & R^{\ot{k}}   &                                    &                                         \\
};
\draw[ ->](m-1-3) edge         node[left ]{$                          $} (m-2-1);
\draw[ ->](m-1-3) edge[dashed] node[left ]{$                          $} (m-2-3);
\draw[ ->](m-1-3) edge         node[left ]{$                          $} (m-2-5);
\draw[ ->](m-2-1) edge         node[right]{$(1\ot\mu)\ot1             $} (m-3-1);
\draw[ ->](m-3-1) edge         node[right]{$                          $} (m-5-1);
\draw[ ->](m-5-1) edge         node[right]{$(1\ot\mu)\ot1             $} (m-6-1);
\draw[ ->](m-6-1) edge         node[right]{$                          $} (m-7-3);
\draw[ ->](m-2-5) edge         node[left ]{$(1\ot\mu)\ot1             $} (m-3-5);
\draw[ ->](m-3-5) edge         node[left ]{$                          $} (m-5-5);
\draw[ ->](m-5-5) edge         node[left ]{$(1\ot\mu)\ot1             $} (m-6-5);
\draw[ ->](m-6-5) edge         node[left ]{$                          $} (m-7-3);
\draw[ ->](m-6-3) edge[dashed] node[right]{$                          $} (m-7-3);
\draw[ ->](m-2-1) edge[dashed] node[right]{$                          $} (m-2-3);
\draw[ ->](m-2-5) edge[dashed] node[right]{$                          $} (m-2-3);
\draw[ ->](m-2-3) edge[sloped,dashed] node[above]{$((1\ot\mu)\ot1)\ot(1\ot1) $} (m-4-2);
\draw[ ->](m-2-3) edge[sloped,dashed] node[above]{$((1\ot1)\ot1)\ot(\mu\ot1) $} (m-4-4);
\draw[ ->](m-4-2) edge[sloped,dashed] node[above]{$((1\ot1)\ot1)\ot(\mu\ot1) $} (m-6-3);
\draw[ ->](m-4-4) edge[sloped,dashed] node[above]{$((1\ot\mu)\ot1)\ot(1\ot1) $} (m-6-3);
\draw[ ->](m-6-1) edge[dashed] node[right]{$                          $} (m-6-3);
\draw[ ->](m-6-5) edge[dashed] node[right]{$                          $} (m-6-3);
\draw[ ->](m-3-1) edge[dashed] node[right]{$                          $} (m-4-2);
\draw[<- ](m-4-2) edge[dashed] node[right]{$                          $} (m-5-1);
\draw[ ->](m-3-5) edge[dashed] node[right]{$                          $} (m-4-4);
\draw[<- ](m-4-4) edge[dashed] node[right]{$                          $} (m-5-5);
\end{tikzpicture}
\end{center}
Again, unlabeled arrows are associativity relations. This diagram obviously commutes by naturality and successive compositions. This proves the claim, and then $(v)=0$. This completes the proof.\mbox{}\hfill\qed\end{mypf}

\subsection{Bimodule objects} \label{sec:bimod}
We take a detour and introduce the Hochschild cochain complex for a module object. In classical Hochschild cohomology the cochain complex has objects of the from $C^{m}(A)=\Hom_{k}(A^{\otimes m},X)$, where $k$ is a field, $A$ a finite dimensional $k$-algebra, and $X=_AX_A$ an $A$-bimodule, i.e.\ the similar setting that we now aim to generalise.

Let $(\fC,\ot,I,\alpha,\lambda,\rho)$ be an $\Ab$-enriched monoidal category and let $(R,\mu,e)$ be a ring object in $\fC$. A \emph{left $R$-module object} is an object $A$ in $\fC$ together with a morphism $\nu:R \ot A \to A$ in $\fC$ called a \emph{left action}, such that the following two diagrams commute 
\newline
\begin{minipage}{.5\textwidth}
\begin{center}
\begin{tikzpicture}
\matrix(m)[matrix of math nodes,row sep=2.4em,column sep=2.8em,text height=1.5ex,text depth=0.25ex]
{
    (R\ot R)\ot A &     &  R \ot(R \ot A)    \\
     R \ot A      &     &  R \ot A           \\
                  & A   &                    \\
};
\draw[ ->,font=\scriptsize](m-1-1) edge         node[above]{$ \alpha_{R,R,A}   $}      (m-1-3);
\draw[ ->,font=\scriptsize](m-1-1) edge         node[right]{$ \mu\ot1_A        $}      (m-2-1);
\draw[ ->,font=\scriptsize](m-2-1) edge         node[above]{$ \nu              $}      (m-3-2);
\draw[ ->,font=\scriptsize](m-1-3) edge         node[left ]{$ 1_R\ot\nu        $}      (m-2-3);
\draw[ ->,font=\scriptsize](m-2-3) edge         node[above]{$ \nu              $}      (m-3-2);
\end{tikzpicture}
\end{center}
\end{minipage}
\begin{minipage}{.5\textwidth}
\begin{center}
\begin{tikzpicture}
\matrix(m)[matrix of math nodes,row sep=2.4em,column sep=2.8em,text height=1.5ex,text depth=0.25ex]
{
    I\ot A &     & R \ot A \\
           & A   &         \\
};
\draw[ ->,font=\scriptsize](m-1-1) edge         node[above]{$e  \ot1_A    $} (m-1-3);
\draw[ ->,font=\scriptsize](m-1-1) edge         node[above]{$\lambda_A    $} (m-2-2);
\draw[ ->,font=\scriptsize](m-1-3) edge         node[above]{$\nu          $} (m-2-2);
\end{tikzpicture}
\end{center}\end{minipage}
We denote this by the pair $(A,\nu)$. Clearly any ring object $R$ is a left module object over itself with left action $\nu=\mu:R\ot R\to R$. We define morphisms of left module objects as follows, let $(A,\nu)$ and $(A',\nu')$ be left module objects. A \emph{morphism of left $R$-module objects} $f:A\to A'$ is a morphism in $\fC$ preserving the action, i.e.\ such that the following diagram commutes 
\begin{center}
\begin{tikzpicture}
\matrix(m)[matrix of math nodes,row sep=2.4em,column sep=2.8em,text height=1.5ex,text depth=0.25ex]
{
    R\ot A &     & R \ot A' \\
    A      &     & A'       \\
};
\draw[ ->,font=\scriptsize](m-1-1) edge         node[above]{$1_R\ot f     $} (m-1-3);
\draw[ ->,font=\scriptsize](m-2-1) edge         node[above]{$f            $} (m-2-3);
\draw[ ->,font=\scriptsize](m-1-1) edge         node[right]{$\nu          $} (m-2-1);
\draw[ ->,font=\scriptsize](m-1-3) edge         node[left ]{$\nu'         $} (m-2-3);
\end{tikzpicture}
\end{center}
Similarly, a \emph{right $R$-module object} $B$ is an object in $\fC$ together with a \emph{right action} $\sigma: B \ot R \to B$ satisfying the required relations given by the following two diagrams 
\newline
\begin{minipage}{.5\textwidth}
\begin{center}
\begin{tikzpicture}
\matrix(m)[matrix of math nodes,row sep=2.4em,column sep=2.8em,text height=1.5ex,text depth=0.25ex]
{
    (B\ot R)\ot R &     &  B \ot(R \ot R)    \\
     B \ot R      &     &  B \ot R           \\
                  & B   &                    \\
};
\draw[ ->,font=\scriptsize](m-1-1) edge         node[above]{$ \alpha_{B,R,R}   $} (m-1-3);
\draw[ ->,font=\scriptsize](m-1-1) edge         node[right]{$ \sigma\ot1_R     $} (m-2-1);
\draw[ ->,font=\scriptsize](m-2-1) edge         node[above]{$ \sigma           $} (m-3-2);
\draw[ ->,font=\scriptsize](m-1-3) edge         node[left ]{$ 1_B\ot\mu        $} (m-2-3);
\draw[ ->,font=\scriptsize](m-2-3) edge         node[above]{$ \sigma           $} (m-3-2);
\end{tikzpicture}
\end{center}
\end{minipage}
\begin{minipage}{.5\textwidth}
\begin{center}
\begin{tikzpicture}
\matrix(m)[matrix of math nodes,row sep=2.4em,column sep=2.8em,text height=1.5ex,text depth=0.25ex]
{
    B\ot I &     & B \ot R \\
           & B   &         \\
};
\draw[ ->,font=\scriptsize](m-1-1) edge         node[above]{$1_B\ot e_R   $} (m-1-3);
\draw[ ->,font=\scriptsize](m-1-1) edge         node[above]{$\rho_B       $} (m-2-2);
\draw[ ->,font=\scriptsize](m-1-3) edge         node[above]{$\sigma       $} (m-2-2);
\end{tikzpicture}
\end{center}
\end{minipage}
Now let $(R,\mu_R,e_R)$ and $(S,\mu_S,e_S)$ be ring objects and let $(X,\nu)$ be a left $S$-module object and $(X,\sigma)$ a right $R$-module object. Then $X$ is said to be an \emph{$(S,R)$-bimodule object} if, in addition, the following diagram commutes 
\begin{center}
\begin{tikzpicture}
\matrix(m)[matrix of math nodes,row sep=2.4em,column sep=2.8em,text height=1.5ex,text depth=0.25ex]
{
    (S \ot X) \ot R &     &  S \ot(X \ot R)    \\
     X \ot R        &     &  S \ot X           \\
                    & X   &                    \\
};
\draw[ ->,font=\scriptsize](m-1-1) edge         node[above]{$ \alpha_{S,X,R}   $} (m-1-3);
\draw[ ->,font=\scriptsize](m-1-1) edge         node[right]{$ \nu\ot1_R        $} (m-2-1);
\draw[ ->,font=\scriptsize](m-2-1) edge         node[above]{$ \sigma           $} (m-3-2);
\draw[ ->,font=\scriptsize](m-1-3) edge         node[left ]{$ 1_S\ot\sigma     $} (m-2-3);
\draw[ ->,font=\scriptsize](m-2-3) edge         node[above]{$ \nu              $} (m-3-2);
\end{tikzpicture}
\end{center}
When $S=R$ we simply say that $(X,\nu,\sigma)$ is an $R$-bimodule object. We remark that any ring object $R$ is an $R$-bimodule over itself where both the left and the right action is given by $\nu=\sigma=\mu:R\ot R\to R$. 

\sloppy We are now going to restate the Hochschild cochain complex for bimodule objects. Let $(\fC,\ot,I,\alpha,\lambda,\rho)$ be an $\Ab$-enriched monoidal category, $(R,\mu,e)$ be a ring object in $\fC$ and let $(X,\nu,\sigma)$ be an $R$-bimodule object. The Hochschild cochain complex $C^{\bullet}(X)=(C^k(X),d^k)_{k\in\mathbb{Z}}$ is defined to be the sequence 
\begin{align*}
\dots \to C^{-1}(X) \xrightarrow{d^{-1}} C^0(X) \xrightarrow{d^0} C^1(X) \xrightarrow {d^1} \cdots 
\end{align*}
that has objects
\begin{align*}
C^k(X)=
\begin{cases}
0                           &\text{for } k<0      \\
\Hom_{\fC}(I,X)             &\text{for } k=0      \\ 
\Hom_{\fC}(R^{\ot k},X) &\text{for } k\geq 1, 
\end{cases}
\end{align*}
The differentials $d^k: C^k(X) \to C^{k+1}(X)$ are defined by
\begin{itemize}
\item $d^k=0$ for $k<0$.
\item For $f\in C^0(X)=\Hom_{\fC}(I,X)$ the differential $d^0:\Hom_{\fC}(I,X)\to\Hom_{\fC}(R,X)$ is defined to be
\begin{align*}
d^0(f)=&\nu \circ (1_R\ot f) \circ \lambda^{-1}_R - \sigma \circ (f\ot 1_R) \circ \rho^{-1}_R
\end{align*}
\item For $k \geq 1$ and $f\in C^k(X) = \Hom_{\fC}(R^{\ot k},X)$ the differentials $d^k:\Hom_{\fC}(R^{\ot k},X)\to\Hom_{\fC}(R^{\ot(k+1)},X)$ are defined to be
\begin{align*}
d^k(f)={          {\nu   } \circ ({1_R \ot f}) \circ \alpha^{0,1}_{k+1} }  
      &+\sum^{k}_{i=1}(-1)^i[{f} \circ ((\alpha^{i-1,1}_k)^{-1} \circ {\mu^{i-1}_{k+1}} \circ \alpha^{i-1,2}_{k+1})] \\
      &+(-1)^{k+1}[{\sigma} \circ {(f \ot 1_R)}] 
\end{align*}\end{itemize}
The sequence $C^{\bullet}(X)=(C^k(X),d^k)_{k\in\mathbb{Z}}$  is indeed a complex i.e.\ $d^{k+1} \circ d^k=0$ for all $k \in \mathbb{Z}$. The proof follows that of Theorem \ref{thm:hcc} by replacing $\mu$ by the actions $\nu$ and $\sigma$ in the places where this makes sense in the diagrams. 

\section{The Hochschild cohomology groups}
We now define Hochschild cohomology as the homology of the cochain complex described in Definition \ref{def:hcc}. Recall that in the classical case of an algebra $A$, the low dimensional groups $\HH^0(A)$, $\HH^1(A)$ and $\HH^2(A)$ can be interpreted in terms of the centre of $A$, the derivations on $A$ and extensions of $A$. For a ring object $R$ we provide similar interpretations for $\HH^0(R)$, $\HH^1(R)$ and $\HH^2(R)$. Throughout this section, we fix an $\Ab$-enriched monoidal category $(\fC,\ot,I,\alpha,\lambda,\rho)$ and a ring object $(R,\mu,e)$ in $\fC$. 

\subsection{The Hochschild cohomology groups}
We define the cohomology groups as follows. 

\begin{mydef} Let $C^{\bullet}(R)=(C^k(R),d^k)_{k\in\mathbb{Z}}$ be the Hochschild cochain complex associated to $R$ (see Definition \ref{def:hcc}). We define the \emph{Hochschild cohomology groups} as the quotient groups 
\begin{align}
\HH^k(R)=\Ker d^k / \IM d^{k-1}, 
\end{align}
i.e.\ the homology of $C^{\bullet}(R)$. \end{mydef}

We remark that these are clearly abelian groups, as quotients of subgroups of abelian groups. Note also remark that the definition of the cohomology groups makes sense for $R$-bimodule objects as in Section \ref{sec:bimod}.  We now discuss some properties of the lower dimensional cohomology groups. 

\subsection{The centre and $\HH^0(R)$}\label{sec:centre}
By definition, 
\begin{align}
\HH^0(R)=\Ker d^0=\{f\in\Hom_{\fC}(I , R) \,|\,\mu \circ (f \ot 1_R) \circ \lambda^{-1}_R-\mu \circ (1_R \ot f) \circ \rho^{-1}_R=0\}. 
\end{align}
With the motivation from the classical case we define the \emph{centre} of $R$ to be $\HH^0(R)$, i.e.\ the collection of $f\in\Hom_{\fC}(I , R)$ such that the following diagram commutes 
\begin{center}
\begin{tikzpicture}
\matrix(m)[matrix of math nodes,row sep=2.6em,column sep=2.6em,text height=1.5ex,text depth=0.25ex]
{   
        & R &          \\
I \ot R &   & R \ot I  \\
R \ot R &   & R \ot R  \\
        & R &          \\
};
\draw[ ->,font=\scriptsize](m-1-2) edge         node[left ]{$\lambda^{-1}_R            $} (m-2-1);
\draw[ ->,font=\scriptsize](m-1-2) edge         node[right]{$\rho^{-1}_R               $} (m-2-3);
\draw[ ->,font=\scriptsize](m-2-1) edge         node[left ]{$f\ot1_R                   $} (m-3-1);
\draw[ ->,font=\scriptsize](m-2-3) edge         node[right]{$1_R\ot f                  $} (m-3-3);
\draw[ ->,font=\scriptsize](m-3-1) edge         node[left ]{$\mu                       $} (m-4-2);
\draw[ ->,font=\scriptsize](m-3-3) edge         node[right]{$\mu                       $} (m-4-2);
\end{tikzpicture}
\end{center}
We denote the centre by $Z(R)$. 

If $R$ happens to be a commutative ring object (to be defined next), then the centre of $R$ consists of all $f\in\Hom_{\fC}(I , R)$, i.e.\ $Z(R)=\Hom_{\fC}(I,R)$. For this purpose we recall the definition of symmetric monoidal categories from \cite[Section VII.7]{mac-98}. A \emph{symmetry} on $\fC$ is a natural isomorphism $\gamma: ?\ot? \rightharpoonup ? \ot? $ with components $\gamma_{A,B}:A \ot B \to B \ot A$ (for every pair $A,B$ of objects in $\fC$) such that the following diagrams, called the \emph{unit coherence}, the \emph{associativity coherence} and the \emph{inverse law}, respectively, commute 
\newline
\begin{minipage}{.25\textwidth}
\begin{tikzpicture}
\matrix(m)[matrix of math nodes,row sep=3em,column sep=1em,text height=1.5ex,text depth=0.25ex]
{ X \ot I &    & X \ot I \\
              & X  &             \\};
\draw[ ->,font=\scriptsize](m-1-1) edge         node[above]{$ \gamma_{X,I}          $} (m-1-3);                       
\draw[ ->,font=\scriptsize](m-1-1) edge         node[left ]{$ \rho_X                $} (m-2-2);                       
\draw[ ->,font=\scriptsize](m-1-3) edge         node[right]{$ \lambda_X             $} (m-2-2);
\end{tikzpicture} 
\end{minipage}
\begin{minipage}{.4\textwidth}
\begin{tikzpicture}
\matrix(m)[matrix of math nodes,row sep=3em,column sep=2em,text height=1.5ex,text depth=0.25ex]
{ 
  (X \ot Y)\ot Z &   & (Y \ot X)\ot Z \\
   X \ot(Y \ot Z)&   &  Y \ot(X \ot Z)\\
  (Y \ot Z)\ot X &   &  Y \ot(Z \ot X)\\
};
\draw[ ->,font=\scriptsize](m-1-1) edge         node[above]{$ \gamma_{X,Y}\ot1_C    $} (m-1-3);                       
\draw[ ->,font=\scriptsize](m-3-1) edge         node[above]{$ \alpha_{Y,Z,X}        $} (m-3-3);                       
\draw[ ->,font=\scriptsize](m-1-1) edge         node[right]{$ \alpha_{X,Y,Z}        $} (m-2-1);                       
\draw[ ->,font=\scriptsize](m-2-1) edge         node[right]{$ \gamma_{X,Y\ot Z}     $} (m-3-1);
\draw[ ->,font=\scriptsize](m-1-3) edge         node[left ]{$ \alpha_{Y,X,Z}        $} (m-2-3);                       
\draw[ ->,font=\scriptsize](m-2-3) edge         node[left ]{$ 1_B\ot\gamma_{X,Z}    $} (m-3-3);
\end{tikzpicture} 
\end{minipage}
\begin{minipage}{.2\textwidth}
\begin{tikzpicture}
\matrix(m)[matrix of math nodes,row sep=3em,column sep=1em,text height=1.5ex,text depth=0.25ex]
{ 
                &  Y \ot X      &               \\
   X \ot Y      &               &  X \ot Y      \\
};
\draw[ ->,font=\scriptsize](m-2-1) edge         node[left ]{$ \gamma_{X,Y}        $} (m-1-2);                       
\draw[ ->,font=\scriptsize](m-1-2) edge         node[right]{$ \gamma_{Y,X}        $} (m-2-3);                       
\draw[ ->,font=\scriptsize](m-2-1) edge         node[below]{$ 1_{X\ot Y}          $} (m-2-3);
\end{tikzpicture} 
\end{minipage}
\newline
A monoidal category with symmetry is called a \emph{symmetric monoidal category} and denoted $(\fC,\ot,I,\alpha,\lambda,\rho,\gamma)$. A ring object $(R,\mu,e)$ in a symmetric monoidal category is said to be a \emph{commutative ring object} if9 $\mu \circ \gamma_{R,R}=\mu$, i.e.\ the following diagram commutes 
\begin{center}
\begin{tikzpicture}
\matrix(m)[matrix of math nodes,row sep=1em,column sep=2.8em,text height=1.5ex,text depth=0.25ex]
{ 
    R\ot R       &            &  R\ot R          \\
                 &  R         &                  \\
};
\draw[ ->,font=\scriptsize](m-1-1) edge         node[above]{$\gamma_{R,R}    $} (m-1-3);
\draw[ ->,font=\scriptsize](m-1-1) edge         node[above]{$\mu             $} (m-2-2);
\draw[ ->,font=\scriptsize](m-1-3) edge         node[above]{$\mu             $} (m-2-2);
\end{tikzpicture}
\end{center}
Observe also that if $R$ and $S$ are commutative ring objects, and $f:R \to S$ a morphism of ring objects, then $f$ preserves the commutativity since $\gamma$ is natural. Hence we do not need any additional structure on $f$. This is similar as for classical rings (i.e.\ ring objects in $\Ab$): the commutativity is preserved just by the axioms of ring homomorphisms. 

In order to prove the claim, that for a commutative ring object $R$ we have $Z(R)=\Hom_{\fC}(I,R)$, we pick an arbitrary $f\in\Hom_{\fC}(I,R)$ and consider the following diagram where dashed arrows with symmetries are added as help lines, 
\begin{center}
\begin{tikzpicture}
\matrix(m)[matrix of math nodes,row sep=2.6em,column sep=2.6em,text height=1.5ex,text depth=0.25ex]
{   
        & R &          \\
I \ot R &   & R \ot I  \\
R \ot R &   & R \ot R  \\
        & R &          \\
};
\draw[ ->,font=\scriptsize](m-1-2) edge         node[left ]{$\lambda^{-1}_R            $} (m-2-1);
\draw[ ->,font=\scriptsize](m-1-2) edge         node[right]{$\rho^{-1}_R               $} (m-2-3);
\draw[ ->,font=\scriptsize](m-2-1) edge         node[left ]{$f\ot1_R                   $} (m-3-1);
\draw[ ->,font=\scriptsize](m-2-3) edge         node[right]{$1_R\ot f                  $} (m-3-3);
\draw[ ->,font=\scriptsize](m-3-1) edge         node[left ]{$\mu                       $} (m-4-2);
\draw[ ->,font=\scriptsize](m-3-3) edge         node[right]{$\mu                       $} (m-4-2);
\draw[ ->,font=\scriptsize](m-2-1) edge[dashed] node[above]{$\gamma_{I,R}              $} (m-2-3);
\draw[ ->,font=\scriptsize](m-3-1) edge[dashed] node[above]{$\gamma_{R,R}              $} (m-3-3);
\end{tikzpicture}
\end{center}
The top triangle commutes by the (inverse of the) unit coherence diagram. The middle square commutes since $\gamma$ is natural. The bottom triangle commutes by the defining relation for commutative ring objects. Since this diagram commutes we conclude that every morphism $f\in\Hom_{\fC}(I , R)$ is in the centre of $R$, whenever $R$ is a commutative ring object. The other direction is obvious, when every $f\in\Hom_{\fC}(I , R)$ is in the centre of $R$ the previous diagram commutes and the bottom triangle tells us that $R$ is commutative.

Now let $\fC$ be not necessarily symmetric, $R$ not necessarily commutative and $(X,\nu,\sigma)$ be an $R$-bimodule object. We define the centre of $X$ as 
\begin{align*} 
Z(X) := \HH^0(R) = \Ker d^0=\{f\in\Hom_{\fC}(I , X) \,|\,\sigma \circ (f \ot 1_R) \circ \lambda^{-1}_R-\nu \circ (1_R \ot f) \circ \rho^{-1}_R=0\}. 
\end{align*}
  
\subsection{Derivations and $\HH^1(R)$}
We define the set of \emph{derivations} on $R$ as  
\begin{align}
\Der(R,R)
:=\{f\in\Hom_{\fC}(R,R) \,|\, f \circ \mu = \mu \circ (1_R\ot f) + \mu \circ (f\ot1_R) \}. 
\end{align}
This is a subgroup of $\Hom_{\fC}(R,R)$. We recall that $d^1(f)
= \mu \circ (1_R\ot f) - f \circ \mu + \mu \circ (f\ot1_R)$. Hence we identify derivations with the kernel of $d^1$, i.e.\ $\Ker d^1=\Der(R,R)$. 

Further we define the set of \emph{inner derivations} as  
\begin{align}
\Der^0(R,R)
:=\{g\in\Hom_{\fC}(R,R)\,|\,f\in\Hom_{\fC}(I,R),\, g=\mu \circ (f\ot 1_R) \circ \lambda^{-1}_R - \mu \circ (1_R\ot f) \circ \rho^{-1}_R\}. 
\end{align}
We immediately observe that $\IM d^0=\Der^0(R,R)$ directly from $d^0(f) = \mu \circ (f\ot 1_R) \circ \lambda^{-1}_R - \mu \circ (1_R\ot f) \circ \rho^{-1}_R$. Moreover, by definition we observe that 
\begin{align} 
\HH^1(R) : = \Ker d^1 / \IM d^0 = \Der(R,R) / \Der^0(R,R)
\end{align}
As in the classical case, we refer to non-zero residue classes in $\HH^1(R)$ as \emph{outer derivations}. 

For bimodule objects we observe that we get the same result when defining 
\begin{align*}
\Der  (R,X)&=\{f\in\Hom_{\fC}(R,X) \,|\, f \circ \mu = \nu \circ (1_R\ot f) + \sigma \circ (f\ot1_R)\} \\
\Der^0(R,X)&=\{g\in\Hom_{\fC}(R,X)\,|\,f\in\Hom_{\fC}(I,X),\, g = \sigma \circ (f\ot 1_R) \circ \lambda^{-1}_R - \nu \circ (1_R\ot f) \circ \rho^{-1}_R \}. 
\end{align*}
So, $\HH^1(X)=\Der(R,X) / \Der^0(R,X)$ as claimed. 

\subsection{Extended algebras and $\HH^2(R)$} 
In this section we assume in addition that $\fC$ has finite coproducts which commute with tensor products, i.e.\ $\fC$ is an additive monoidal category. Let $f\in\Hom_{\fC}(R^{\ot2},R)$. When $f\in\Ker d^2$ we are going to construct a ring object $R \ltimes_f R$ in $\fC$ referred to as the \emph{extended ring object of $R$ along $f$}. As object $R \ltimes_f R$ is the coproduct $R \amalg R$ which is equipped with a multiplication rule and a multiplicative unit. 

We recall that morphisms between finite coproducts in additive categories can be represented by matrix expressions where the entries are morphisms between the respective objects (see \cite[Section VIII.2]{mac-98}). In this section we simplify the following notation, when not otherwise stated let $1_R=:1$, $\mu_R=:\mu$, $\alpha_{R,R,R}=:\alpha$, $e_R=:e$ and $\alpha_{R\ltimes_f R,R\ltimes_f R,R\ltimes_f R}=:\hat{\alpha}$. As usual let $0$ denote the unique map that factors uniquely through the unique zero object of $\fC$. 

The distributive law from \cite[p. 172]{mac-98} (in the proof of Theorem 2) will be useful. From this we have that $A \ot (B_1 \amalg B_2) \cong (A \ot B_1) \amalg (A \ot B_2)$. Hence 
\begin{align}
(R \amalg R) \ot (R \amalg R) \cong (R \ot R) \amalg (R \ot R) \amalg (R \ot R) \amalg (R \ot R). 
\end{align}

We define the multiplication rule $\mu_{R \ltimes_f R}:(R \amalg R) \ot (R \amalg R)\to R \amalg R$ on $R \ltimes_f R$ as the composition 
\begin{align*}
\mu_{R \ltimes_f R}:
(R \amalg R) \ot (R \amalg R) 
\xrightarrow[\sim]{\begin{pmatrix} (1,0)\ot(1,0) \\ (1,0)\ot(0,1) \\ (0,1)\ot(1,0) \\ (0,1)\ot(0,1) \end{pmatrix}} 
(R \ot R) \amalg (R \ot R) \amalg (R \ot R) \amalg (R \ot R) 
\\
\xrightarrow[]{\begin{pmatrix}\mu & 0 & 0 & 0 \\ f & \mu & \mu & 0 \end{pmatrix}}
R \amalg R
\end{align*}
which written together is 
\begin{align}
\mu_{R\ltimes_fR}=
\begin{pmatrix}
\mu \circ ((1,0)\ot(1,0))\\
f \circ ((1,0)\ot(1,0)) + \mu \circ ((1,0)\ot(0,1)) + \mu \circ ((0,1)\ot(1,0))
\end{pmatrix}:\nonumber\\(R \amalg R) \ot (R \amalg R) \to  R \amalg R 
\end{align}
We claim that the multiplicative unit $e_{R\ltimes_f R}:I \to R\ltimes_f R$ is defined as the composition, 
\begin{align*}
e_{R\ltimes_f R}:I
\xrightarrow[]{\begin{pmatrix}1_I\\1_I\end{pmatrix}}
I \amalg I 
\xrightarrow[]{\begin{pmatrix}1_I & 0 \\ 0 & \lambda_I^{-1}\end{pmatrix}}
I \amalg (I\ot I)
\xrightarrow[]{\begin{pmatrix}e_R & 0 \\ 0 & e_R\ot e_R\end{pmatrix}}
R \amalg (R \ot R)
\\\xrightarrow[]{\begin{pmatrix}1_R & 0 \\ 0 & -f\end{pmatrix}}
R \amalg R
\end{align*}
Recall that $\lambda_I=\rho_I$, so we can replace them with one another appropriate to the setting we are studying. Composed together this gives, 
\begin{align}
e_{R\ltimes_f R}=
\begin{pmatrix}
e_R \\
-f \circ (e_R \ot e_R) \circ \lambda_I^{-1}
\end{pmatrix}: I \to R \amalg R 
\end{align}

\begin{myprop}
Let $(\fC,\ot,I,\alpha,\lambda,\rho)$ be an additive monoidal category, and $(R,\mu,e)$ a ring object in $\fC$. Then the triple $(R \ltimes_f R , \mu_{R \ltimes_f R} , e_{R \ltimes_f R})$ is a ring object in $\fC$. \end{myprop}

\begin{mypf}
We have to show that the multiplication rule and the multiplicative unit satisfy the required associativity and unitary relations. For the associativity relation we want to show that the following diagram commutes 
\begin{center}
\begin{tikzpicture}
\matrix(m)[matrix of math nodes,row sep=2.6em,column sep=2.8em,text height=1.5ex,text depth=0.25ex]
{(R\ltimes_f R \ot R\ltimes_f R)\ot R\ltimes_f R  &               & R\ltimes_f R \ot(R\ltimes_f R \ot R\ltimes_f R)  \\
  R\ltimes_f R           \ot R\ltimes_f R         &               & R\ltimes_f R \ot R\ltimes_f R                    \\
                                                  &  R\ltimes_f R &                                                  \\};
\draw[ ->,font=\scriptsize](m-1-1) edge         node[above]{$\alpha_{R\ltimes_f R,R\ltimes_f R,R\ltimes_f R}          $} (m-1-3);
\draw[ ->,font=\scriptsize](m-1-1) edge         node[left ]{$\mu_{R\ltimes_f R} \ot 1_{R\ltimes_f R}                  $} (m-2-1);
\draw[ ->,font=\scriptsize](m-2-1) edge         node[below]{$\mu_{R\ltimes_f R}                                       $} (m-3-2);

\draw[ ->,font=\scriptsize](m-1-3) edge         node[right]{$1_{R\ltimes_f R}   \ot \mu_{R\ltimes_f R}                $} (m-2-3);
\draw[ ->,font=\scriptsize](m-2-3) edge         node[below]{$\mu_{R\ltimes_f R}                                       $} (m-3-2);
\end{tikzpicture}
\end{center}
Composing the ``left'' side we get, 
\begin{align*}
&
\mu_{R\ltimes_f R} \circ (\mu_{R\ltimes_f R} \ot 1_{R\ltimes_f R})
\\&=
\begin{pmatrix}
\mu & 0  & 0  & 0 \\
f   &\mu &\mu & 0 
\end{pmatrix}
\begin{pmatrix}
(1,0)\ot(1,0) \\
(1,0)\ot(0,1) \\
(0,1)\ot(1,0) \\
(0,1)\ot(0,1)  
\end{pmatrix}
\begin{pmatrix}
\begin{pmatrix}
\mu & 0  & 0  & 0 \\
f   &\mu &\mu & 0 
\end{pmatrix}
\begin{pmatrix}
(1,0)\ot(1,0) \\
(1,0)\ot(0,1) \\
(0,1)\ot(1,0) \\
(0,1)\ot(0,1)  
\end{pmatrix}
\ot
\begin{pmatrix}
1 & 0 \\
0 & 1 
\end{pmatrix}
\end{pmatrix}
\\&= 
\begin{pmatrix}
\mu \circ ( [\mu \circ ((1,0)\ot(1,0))] \ot (1,0))\\
\xi_1
\end{pmatrix}
=:\mathpzc{L}
\end{align*}
where 
\begin{align*}
\xi_1=&
\overbrace{f  \circ ([\mu \circ ((1,0)\ot(1,0))]\ot(1,0))}^{(1)} + 
\overbrace{\mu\circ ([\mu \circ ((1,0)\ot(1,0))]\ot(0,1))}^{(2)} + \\& 
\overbrace{\mu\circ ([f   \circ ((1,0)\ot(1,0))]\ot(1,0))}^{(3)} + 
\overbrace{\mu\circ ([\mu \circ ((1,0)\ot(0,1))]\ot(1,0))}^{(4)} + \\&
\overbrace{\mu\circ ([\mu \circ ((0,1)\ot(1,0))]\ot(1,0))}^{(5)}
\end{align*}
For the other path, namely compositions of arrows on the upper and right direction, we get the following
\begin{align*}
&
\mu_{R\ltimes_f R} \circ (\mu_{R\ltimes_f R} \ot 1_{R\ltimes_f R}) \circ \alpha_{R\ltimes_fR,R\ltimes_fR,R\ltimes_fR}=\mu_{R\ltimes_f R} \circ (\mu_{R\ltimes_f R} \ot 1_{R\ltimes_f R}) \circ \hat{\alpha}
\\&=
\begin{pmatrix}
\mu & 0  & 0  & 0 \\
f   &\mu &\mu & 0 
\end{pmatrix}
\begin{pmatrix}
(1,0)\ot(1,0) \\
(1,0)\ot(0,1) \\
(0,1)\ot(1,0) \\
(0,1)\ot(0,1)  
\end{pmatrix}
\begin{pmatrix}
\begin{pmatrix}
1 & 0 \\
0 & 1 
\end{pmatrix}
\ot
\begin{pmatrix}
\mu & 0  & 0  & 0 \\
f   &\mu &\mu & 0 
\end{pmatrix}
\begin{pmatrix}
(1,0)\ot(1,0) \\
(1,0)\ot(0,1) \\
(0,1)\ot(1,0) \\
(0,1)\ot(0,1)  
\end{pmatrix}
\end{pmatrix}\hat{\alpha}
\\&= 
\begin{pmatrix}
\mu \circ ((1,0) \ot [\mu \circ ((1,0)\ot(1,0))                                         ])\\
\xi_2
\end{pmatrix}\hat{\alpha}
=:\mathpzc{R}
\end{align*}
where 
\begin{align*}
\xi_2=&
\overbrace{(f  \circ((1,0)\ot [\mu \circ ((1,0)\ot(1,0))]))\hat{\alpha}}^{(1')} + \\&
\overbrace{(\mu\circ((1,0)\ot [f   \circ ((1,0)\ot(1,0))]))\hat{\alpha}}^{(2')} + 
\overbrace{(\mu\circ((1,0)\ot [\mu \circ ((1,0)\ot(0,1))]))\hat{\alpha}}^{(3')} + \\&
\overbrace{(\mu\circ((1,0)\ot [\mu \circ ((0,1)\ot(1,0))]))\hat{\alpha}}^{(4')} + 
\overbrace{(\mu\circ((0,1)\ot [\mu \circ ((1,0)\ot(1,0))]))\hat{\alpha}}^{(5')}
\end{align*}
First we check if the top entries in $\mathpzc{L}$ and $\mathpzc{R}$ coincide. Consider the diagram 
\begin{center}
\begin{tikzpicture}
\matrix(m)[matrix of math nodes,row sep=2.6em,column sep=2.8em,text height=1.5ex,text depth=0.25ex]
{(R\amalg R \ot R\amalg R)\ot R\amalg R  &               & R\amalg R \ot(R\amalg R \ot R\amalg R)  \\
 (R         \ot R        )\ot R          &               & R         \ot(R         \ot R        )  \\
  R                       \ot R          &               & R         \ot R                         \\
                                         &  R            &                                         \\};
\draw[ ->,font=\scriptsize](m-1-1) edge         node[above]{$\alpha_{R\amalg R,R\amalg R,R\amalg R}          $} (m-1-3);
\draw[ ->,font=\scriptsize](m-1-1) edge         node[left ]{$((1,0)\ot (1,0))\ot(1,0)                        $} (m-2-1);
\draw[ ->,font=\scriptsize](m-1-3) edge         node[left ]{$ (1,0)\ot((1,0) \ot(1,0))                       $} (m-2-3);
\draw[ ->,font=\scriptsize](m-2-1) edge         node[above]{$\alpha_{R,R,R}                                  $} (m-2-3);

\draw[ ->,font=\scriptsize](m-2-1) edge         node[left ]{$\mu\ot1                                         $} (m-3-1);
\draw[ ->,font=\scriptsize](m-3-1) edge         node[below]{$\mu                                             $} (m-4-2);
\draw[ ->,font=\scriptsize](m-2-3) edge         node[left ]{$1\ot\mu                                         $} (m-3-3);
\draw[ ->,font=\scriptsize](m-3-3) edge         node[below]{$\mu                                             $} (m-4-2);
\end{tikzpicture}
\end{center}
which commutes: the top by naturality of $\alpha$ and the bottom by the associativity relation of $\mu$. Then we have that  
\begin{align*}
\mu\circ([\mu \circ((1,0)\ot(1,0))]\ot(1,0))=\mu\circ((1,0)\ot [\mu \circ((1,0)\ot(1,0))])\circ \alpha_{R\ltimes_f R,R\ltimes_f R,R\ltimes_f R},
\end{align*}
which proves that the top entries in $\mathpzc{L}$ and $\mathpzc{R}$ coincide. 
Next we prove that the bottom entries in $\mathpzc{L}$ and $\mathpzc{R}$ coincide, i.e.\ $\xi_1=\xi_2$. 
We use that $\alpha$ is natural to identify $(2)=(3')$, $(4)=(4')$ and $(5)=(5')$. For the remaining terms we determine   
\begin{align*}
&
d^2(f) \circ (((1,0)\ot(1,0))\ot(1,0))
\\&=
(\mu \circ (1 \ot f) \circ \alpha_{R,R,R} - f \circ (\mu \ot 1) + f \circ (1 \ot \mu) \circ \alpha_{R,R,R} - \mu \circ (f \ot 1)) \circ (((1,0)\ot(1,0))\ot(1,0))
\\&=
\mu \circ ((1,0) \ot [f \circ ((1,0)\ot(1,0))]) \circ \hat{\alpha} - f \circ ([\mu \circ ((1,0)\ot(1,0))] \ot (1,0)) + \\&\hspace{3cm} f \circ ((1,0) \ot [\mu \circ ((1,0)\ot(1,0))]) \circ \hat{\alpha} - \mu \circ ([f \circ ((1,0)\ot(1,0))] \ot (1,0))=0 
\end{align*}
which vanishes since $f\in\Ker d^2$, and $\alpha_{R,R,R}\circ(((1,0)\ot(1,0))\ot(1,0))=((1,0)\ot((1,0)\ot(1,0)))\circ\hat{\alpha}$ (since $\alpha$ is natural, i.e.\ the top of the previous diagram commutes). Then 
\begin{align*}
(1')+(2')=&
\overbrace{(f  \circ ((1,0)\ot [\mu \circ((1,0)\ot(1,0))]))\circ\hat{\alpha}}^{(1')} + 
\overbrace{(\mu\circ ((1,0)\ot [f   \circ((1,0)\ot(1,0))]))\circ\hat{\alpha}}^{(2')}  
\\=&
\overbrace{ f  \circ ([\mu \circ((1,0)\ot(1,0))]\ot(1,0))}^{(1)} + 
\overbrace{ \mu\circ ([f   \circ((1,0)\ot(1,0))]\ot(1,0))}^{(3)}
=(1)+(3)
\end{align*}
This finial argument shows that $\xi_1=\xi_2$, and we conclude that $\mathpzc{L}=\mathpzc{R}$ and the multiplication rule $\mu_{R\ltimes_fR}$ is associative. 

Next we show that the unitary relations hold, starting with the left unitary relation. We want to show that
\begin{center}
\begin{tikzpicture}
\matrix(m)[matrix of math nodes,row sep=2.6em,column sep=5.8em,text height=1.5ex,text depth=0.25ex]
{
  R\ltimes_f R & I \ot (R\ltimes_f R) & (R\ltimes_f R) \ot (R\ltimes_f R) \\
               &                      &  R\ltimes_f R                     \\
};
\draw[ ->,font=\scriptsize](m-1-1) edge         node[above]{$\lambda_{R\ltimes_fR}^{-1}               $} (m-1-2);
\draw[ ->,font=\scriptsize](m-1-2) edge         node[above]{$e_{R\ltimes_f R} \ot 1_{R\ltimes_f R}    $} (m-1-3);
\draw[ ->,font=\scriptsize](m-1-3) edge         node[right]{$\mu_{R\ltimes_f R}                       $} (m-2-3);
\draw[ ->,font=\scriptsize](m-1-1) edge         node[below]{$1_{R\ltimes_f R}                         $} (m-2-3);
\end{tikzpicture}
\end{center}
commutes, i.e.\ we want to show that  
$\mu_{R\ltimes_f R}\circ(e_{R\ltimes_f R} \ot 1_{R\ltimes_f R})\circ\lambda_{R\ltimes_f R}^{-1}=1_{R\ltimes_fR}$. 
We compose and calculate  
\begin{align*}
&
\mu_{R\ltimes_f R}\circ(e_{R\ltimes_f R} \ot 1_{R\ltimes_f R})\circ\lambda_{R\ltimes_f R}^{-1}
\\&=
\begin{pmatrix}
\mu & 0   & 0   & 0   \\
f   & \mu & \mu & 0
\end{pmatrix}
\begin{pmatrix}
(1,0)\ot(1,0)\\
(1,0)\ot(0,1)\\
(0,1)\ot(1,0)\\
(0,1)\ot(0,1)  
\end{pmatrix}
\begin{pmatrix}
\begin{pmatrix}
e \\
-f \circ (e \ot e)\circ \lambda_{I}^{-1}
\end{pmatrix}
\ot
\begin{pmatrix}
1 & 0 \\ 
0 & 1 
\end{pmatrix}
\end{pmatrix}
\lambda_{R\ltimes_f R}^{-1}
\\&=
\begin{pmatrix}
\mu\circ( e                                            \ot1)\circ(1_I\ot(1,0))\\
f  \circ[ e                                            \ot(1,0)] + 
\mu\circ[ e                                            \ot1]\circ(1_I\ot(0,1)) -  
\mu\circ[(   f \circ (e \ot e) \circ  \lambda_{I}^{-1})\ot(1,0)] 
\end{pmatrix}
\lambda_{R\ltimes_f R}^{-1}
\\&=
\begin{pmatrix}
\lambda\circ(1_I\ot(1,0))\\
\underbrace{f  \circ [ e                          \ot(1,0)]}_{(i)} + 
            \lambda\circ(1_I\ot(0,1)) -  
\underbrace{\mu\circ[(   f \circ(e \ot e)\circ \lambda_{I}^{-1})\ot(1,0)]}_{(ii)} 
\end{pmatrix}
\lambda_{R\ltimes_f R}^{-1}
\end{align*}
where we have used that $\mu\circ(e\ot1)=\lambda$. We claim that $(i)-(ii)=f\circ[e\ot(1,0)] - \mu\circ[(f \circ(e \ot e)\circ \lambda_{I}^{-1})\ot(1,0)] = 0$, in fact this is $d^2(f) \circ ((e \ot (e \circ \lambda_I^{-1})) \ot (1,0))$ which clearly vanishes, since $f\in\Ker d^2$. To see this we calculate 
\begin{align*}
&d^2(f) \circ ((e \ot e) \circ \lambda_I^{-1}) \ot (1,0))=d^2(f)\circ((e \ot e) \ot 1)\circ(\lambda_I^{-1} \ot (1,0))
\\&=
(\mu \circ [1 \ot f] \circ \alpha - f \circ [\mu \ot 1] + f \circ [1 \ot \mu] \circ \alpha - \mu \circ [f \ot 1] )\circ((e \ot e) \ot 1)\circ(\lambda_I^{-1} \ot (1,0))
\\&= 
\left(
\overbrace{\mu \circ [e \ot f \circ (e \ot 1)] \circ \alpha_{I,I,R}  }^{(1)} - 
\overbrace{f   \circ [\mu \circ (e\ot e) \ot 1]                }^{(2)} + 
\overbrace{f   \circ [e \ot \mu \circ (e \ot 1)] \circ \alpha_{I,I,R}}^{(3)} -\right.\\ &\hspace{3cm}\left.
\overbrace{\mu \circ [ f \circ (e \ot e)\ot 1]               }^{(4)}   
\right) \circ (\lambda_I^{-1} \ot (1,0))
=0. 
\end{align*}
We rewrite slightly, 
\begin{align*}
&(1)=\mu \circ [e \ot f \circ (e \ot1)] \circ \alpha_{I,I,R} = \mu \circ (e\ot1) \circ (1\ot f \circ (e \ot 1))\circ\alpha_{I,I,R} \\&\hspace{8cm}= \lambda \circ (1_I\ot f \circ (e\ot1))\circ\alpha_{I,I,R}\\  
&(2)=f \circ [\mu \circ (e\ot e) \ot 1]              =f\circ[\mu\circ(e\ot1)\circ(1\ot e)\ot1]            =f\circ[\lambda\circ(1\ot e)\ot1]             \\
&(3)=f \circ [e \ot \mu\circ (e \ot1)]\circ\alpha_{I,I,R}=f \circ [e \ot \lambda]\circ \alpha_{I,I,R}
\end{align*}
Later we will identify $(i)-(ii)=(3)-(4)$. With this identification assumed to be true, and if we in addition can show that 
\begin{align*}
(1)-(2)=\lambda \circ (1\ot f \circ (e\ot1))\circ \alpha_{I,I,R}-f \circ [\lambda \circ (1\ot e)\ot1]=0
\end{align*}
we are done. We approach this latter by considering the following diagram 
\begin{center}
\begin{tikzpicture}
\matrix(m)[matrix of math nodes,row sep=2.6em,column sep=2.8em,text height=1.5ex,text depth=0.25ex]
{
 (I \ot  I) \ot R      &  I \ot (I  \ot R)&              \\
 (I \ot  R) \ot R      &  I \ot (R  \ot R)& I \ot R      \\
  R \ot  R             &  R \ot  R        & R            \\
};
\draw[ ->,font=\scriptsize](m-1-1) edge         node[above]{$\alpha_{I,I,R}         $} (m-1-2);
\draw[ ->,font=\scriptsize](m-1-2) edge         node[right]{$ 1_I\ot(e \ot 1)       $} (m-2-2);
\draw[ ->,font=\scriptsize](m-1-1) edge         node[left ]{$(1_I\ot e)\ot 1        $} (m-2-1);
\draw[ ->,font=\scriptsize](m-2-1) edge         node[below]{$\alpha_{I,R,R}         $} (m-2-2);
\draw[ ->,font=\scriptsize](m-2-1) edge         node[left ]{$\lambda\ot1               $} (m-3-1);
\draw[ ->,font=\scriptsize](m-3-1) edge         node[below]{$1_{R\ot R}             $} (m-3-2);
\draw[ ->,font=\scriptsize](m-2-2) edge         node[right]{$\lambda_{R\ot  R}         $} (m-3-2);
\draw[ ->,font=\scriptsize](m-2-2) edge         node[above]{$1_I\ot f               $} (m-2-3);
\draw[ ->,font=\scriptsize](m-3-2) edge         node[below]{$f                      $} (m-3-3);
\draw[ ->,font=\scriptsize](m-2-3) edge         node[right]{$\lambda                   $} (m-3-3);
\end{tikzpicture}
\end{center}
where the outer/top zigzag is $(1)$ and the left/bottom edges are $(2)$. The top square commutes since $\alpha$ is a natural transformation, the bottom left square commutes by the coherence theorem, the bottom right square commutes since $\lambda$ is an natural transformation. So $(1)-(2)=0$. Hence we are left with 
\begin{align*}
d^2(f)\circ((e \ot e) \ot 1)\circ(\lambda_I^{-1} \ot (1,0))
&=((3)-(4))\circ(\lambda_I^{-1} \ot (1,0))=0
\end{align*}
We identify 
\begin{align*} 
(4)\circ(\lambda_I^{-1} \ot (1,0))=\mu \circ [f \circ (e \ot e)\ot 1] \circ(\lambda_I^{-1} \ot (1,0))=\mu \circ [(f \circ (e \ot e) \circ \lambda_I^{-1})\ot (1,0)]=(ii). 
\end{align*}
Further we rewrite 
\begin{align*}
(3)\circ(\lambda_I^{-1} \ot (1,0))&=f \circ [e \ot \lambda]\circ\alpha_{I,I,R}\circ(\lambda_I^{-1} \ot (1,0))\\
                                  &=f \circ (e\ot1)\circ(1 \ot \lambda)\circ\alpha_{I,I,R}\circ(\lambda_I^{-1} \ot 1)\circ(1_I\ot(1,0)). 
\end{align*}
Consider the following diagram 
\begin{center}
\begin{tikzpicture}
\matrix(m)[matrix of math nodes,row sep=2.6em,column sep=2.8em,text height=1.5ex,text depth=0.25ex]
{
 (I \ot  I) \ot R      &  I \ot (I  \ot R)& I \ot (I  \ot R)             \\
  I \ot  R             &  I \ot  R        & I \ot R                      \\
};
\draw[ ->,font=\scriptsize](m-1-1) edge         node[above]{$\alpha_{I,I,R}           $} (m-1-2);
\draw[ ->,font=\scriptsize](m-1-2) edge         node[above]{$1_{I \ot (I  \ot R)}     $} (m-1-3);
\draw[ ->,font=\scriptsize](m-2-1) edge         node[above]{$1_{I \ot  R }            $} (m-2-2);
\draw[ ->,font=\scriptsize](m-2-2) edge         node[above]{$1_{I \ot  R }            $} (m-2-3);
\draw[ ->,font=\scriptsize](m-1-1) edge         node[left ]{$\lambda_I\ot1               $} (m-2-1);
\draw[ ->,font=\scriptsize](m-1-2) edge         node[right]{$\lambda_{I\ot R}            $} (m-2-2);
\draw[ ->,font=\scriptsize](m-1-3) edge         node[right]{$1_I\ot\lambda               $} (m-2-3);
\end{tikzpicture}
\end{center}
which commutes by the coherence theorem. We obtain $(1 \ot \lambda)\circ\alpha_{I,I,R}\circ(\lambda_I^{-1} \ot 1)=1_{I \ot R}$, which gives 
\begin{align*}
f \circ (e\ot1)\circ(1 \ot \lambda)\circ\alpha_{I,I,R}\circ(\lambda_I^{-1} \ot 1)\circ(1_I\ot(1,0))=f\circ(e\ot(1,0))=(i). 
\end{align*}
Hence 
\begin{align*}
0=d^2(f)\circ((e \ot e) \ot 1)\circ(\lambda_I^{-1} \ot (1,0))&=f\circ[e\ot(1,0)] - \mu\circ[(f \circ (e \ot e) \circ \lambda_{I}^{-1})\ot(1,0)]\\&=(i)-(ii). 
\end{align*}
as claimed. So far in the left unitary law we are left with 
\begin{align}
\mu_{R\ltimes_f R}\circ(e_{R\ltimes_f R} \ot 1_{R\ltimes_f R})\circ\lambda_{R\ltimes_f R}^{-1}
=
\begin{pmatrix}
\lambda(1_I\ot(1,0))\\
\lambda(1_I\ot(0,1))
\end{pmatrix}
\lambda_{R\ltimes_f R}^{-1}
\end{align}
Now consider the following two diagrams 
\newline
\begin{minipage}{.5\textwidth}
\begin{center}
\begin{tikzpicture}
\matrix(m)[matrix of math nodes,row sep=2.6em,column sep=2.8em,text height=1.5ex,text depth=0.25ex]
{
  I \ot (R\ltimes_f R) & R\ltimes_f R \\
  I \ot  R             & R            \\
};
\draw[ ->,font=\scriptsize](m-1-1) edge         node[above]{$\lambda_{R\ltimes_fR}     $} (m-1-2);
\draw[ ->,font=\scriptsize](m-1-2) edge         node[left ]{$(1,0)                  $} (m-2-2);
\draw[ ->,font=\scriptsize](m-1-1) edge         node[right]{$1_I\ot(1,0)            $} (m-2-1);
\draw[ ->,font=\scriptsize](m-2-1) edge         node[below]{$\lambda                   $} (m-2-2);
\end{tikzpicture}
\end{center}
\end{minipage}
\begin{minipage}{.5\textwidth}
\begin{center}
\begin{tikzpicture}
\matrix(m)[matrix of math nodes,row sep=2.6em,column sep=2.8em,text height=1.5ex,text depth=0.25ex]
{
  I \ot (R\ltimes_f R) & R\ltimes_f R \\
  I \ot  R             & R            \\
};
\draw[ ->,font=\scriptsize](m-1-1) edge         node[above]{$\lambda_{R\ltimes_fR}     $} (m-1-2);
\draw[ ->,font=\scriptsize](m-1-2) edge         node[left ]{$(0,1)                  $} (m-2-2);
\draw[ ->,font=\scriptsize](m-1-1) edge         node[right]{$1_I\ot(0,1)            $} (m-2-1);
\draw[ ->,font=\scriptsize](m-2-1) edge         node[below]{$\lambda                   $} (m-2-2);
\end{tikzpicture}
\end{center}
\end{minipage}
which commute since $\lambda$ is natural. Hence we have that $\lambda\circ(1_I\ot(1,0))=(1,0)\circ\lambda_{R\ltimes_f R}$ and $\lambda\circ(1_I\ot(0,1))=(0,1)\circ\lambda_{R\ltimes_f R}$. Insert this back in the original expression and get, 
\begin{align*}
&
\mu_{R\ltimes_f R}\circ(e_{R\ltimes_f R} \ot 1_{R\ltimes_f R})\circ\lambda_{R\ltimes_f R}^{-1}
\\=&
\begin{pmatrix}
\lambda \circ (1_I\ot(1,0))\\
\lambda \circ (1_I\ot(0,1))
\end{pmatrix}
\lambda_{R\ltimes_f R}^{-1}
=
\begin{pmatrix}
(1,0) \circ \lambda_{R\ltimes_f R}\\
(0,1) \circ \lambda_{R\ltimes_f R}
\end{pmatrix}
\lambda_{R\ltimes_f R}^{-1}
=
\begin{pmatrix}
1&0\\
0&1
\end{pmatrix}
=
1_{R\ltimes_f R}
\end{align*}
which proves the left unitary law. 

For the right unitary law we have to show that
\begin{center}
\begin{tikzpicture}
\matrix(m)[matrix of math nodes,row sep=2.6em,column sep=5.8em,text height=1.5ex,text depth=0.25ex]
{
  R\ltimes_f R & (R\ltimes_f R) \ot I & (R\ltimes_f R) \ot (R\ltimes_f R) \\
               &                      &  R\ltimes_f R                     \\
};
\draw[ ->,font=\scriptsize](m-1-1) edge         node[above]{$\rho_{R\ltimes_fR}^{-1}                  $} (m-1-2);
\draw[ ->,font=\scriptsize](m-1-2) edge         node[above]{$1_{R\ltimes_f R} \ot e_{R\ltimes_f R}    $} (m-1-3);
\draw[ ->,font=\scriptsize](m-1-3) edge         node[right]{$\mu_{R\ltimes_f R}                       $} (m-2-3);
\draw[ ->,font=\scriptsize](m-1-1) edge         node[below]{$1_{R\ltimes_f R}                         $} (m-2-3);
\end{tikzpicture}
\end{center}
commutes. So we have to show that
$\mu_{R\ltimes_f R}(1_{R\ltimes_f R} \ot e_{R\ltimes_f R} )\rho_{R\ltimes_fR}^{-1}=1_{R\ltimes_f R}$. 
We identify $\rho_I=\lambda_I$. Then we calculate 
\begin{align*}
&
\mu_{R\ltimes_f R} \circ (1_{R\ltimes_f R} \ot e_{R\ltimes_f R}) \circ \rho_{R\ltimes_fR}^{-1}
\\&=
\begin{pmatrix}
\mu & 0   & 0   & 0   \\
f   & \mu & \mu & 0
\end{pmatrix}
\begin{pmatrix}
(1,0)\ot(1,0)\\
(1,0)\ot(0,1)\\
(0,1)\ot(1,0)\\
(0,1)\ot(0,1)  
\end{pmatrix}
\begin{pmatrix}
\begin{pmatrix}
1 & 0 \\ 
0 & 1 
\end{pmatrix}
\ot
\begin{pmatrix}
e \\
-f \circ (e \ot e) \circ \rho_{I}^{-1}
\end{pmatrix}
\end{pmatrix}
\rho_{R\ltimes_f R}^{-1}
\\&=
\begin{pmatrix}
\mu \circ ((1,0)\ot e) \\
f \circ ((1,0)\ot e) - \mu \circ ((1,0) \ot (f \circ (e \ot e) \circ \rho_{I}^{-1})) + \mu \circ ((0,1)\ot e) 
\end{pmatrix}\rho_{R\ltimes_f R}^{-1}
\end{align*}
Similarly as for the left unitary law, we will prove that $f \circ ((1,0)\ot e) - \mu ((1,0) \ot f \circ (e \ot e) \circ \rho_{I}^{-1})=0$ from the property that $f\in\Ker d^2$. We have 
\begin{align*}
&
d^2(f)\circ((1 \ot e) \ot e)\circ\alpha_{R,I,I}^{-1}((1,0) \ot \rho_I^{-1})
\\&=
\big(\mu \circ [1 \ot f] \circ \alpha - f \circ [\mu \ot 1] + f \circ [1 \ot \mu] \circ \alpha - \mu \circ [f \ot 1] \big)\circ ((1 \ot e) \ot e) \circ \alpha_{R,I,I}^{-1}\circ ((1,0) \ot \rho_I^{-1})
\\&
=\big(\mu \circ [1 \ot (f \circ (e \ot e))] \circ \alpha_{R,I,I} - f \circ [\rho \ot e]\\&\hspace{1cm} + f \circ (1 \ot \rho) \circ (1\ot(e \ot 1_I)) \circ \alpha_{R,I,I} - \rho \circ (f\ot1)\circ ((1\ot e)\ot1_I)\big) \circ \alpha_{R,I,I}^{-1} \circ ((1,0) \ot \rho_I^{-1})
\end{align*}
We identify the two last terms in the expression above with each other by the following diagram 
\begin{center}
\begin{tikzpicture}
\matrix(m)[matrix of math nodes,row sep=2.6em,column sep=2.8em,text height=1.5ex,text depth=0.25ex]
{
    (R\ot I)\ot I                     &   R \ot (I \ot I)                 \\
    (R\ot R)\ot I                     &   R \ot (R \ot I)                 \\
    (R\ot R)\ot I                     &          R \ot R                  \\
     R      \ot I                     &   R                               \\
};
\draw[ ->,font=\scriptsize](m-1-1) edge         node[above]{$\alpha_{R,I,I}                           $} (m-1-2);
\draw[ ->,font=\scriptsize](m-2-1) edge         node[above]{$\alpha_{R,R,I}                           $} (m-2-2);
\draw[ ->,font=\scriptsize](m-3-1) edge         node[above]{$\rho_{R\ot R}                            $} (m-3-2);
\draw[ ->,font=\scriptsize](m-4-1) edge         node[above]{$\rho                                     $} (m-4-2);
\draw[ ->,font=\scriptsize](m-1-1) edge         node[left ]{$(1\ot e)\ot1_I                           $} (m-2-1);
\draw[ ->,font=\scriptsize](m-2-1) edge         node[left ]{$1_{(R\ot R)\ot I}                        $} (m-3-1);
\draw[ ->,font=\scriptsize](m-3-1) edge         node[left ]{$f \ot 1_I                                $} (m-4-1);
\draw[ ->,font=\scriptsize](m-1-2) edge         node[right]{$ 1\ot(e\ot1_I)                           $} (m-2-2);
\draw[ ->,font=\scriptsize](m-2-2) edge         node[right]{$1\ot \rho                                $} (m-3-2);
\draw[ ->,font=\scriptsize](m-3-2) edge         node[right]{$f                                        $} (m-4-2);
\end{tikzpicture}
\end{center}
where the top square commutes by naturality of $\alpha$, the middle square commutes by the coherence theorem and the bottom square commutes by the naturality of $\rho$. Since the top and right composition is the summand $f \circ (1 \ot \rho) \circ (1\ot(e \ot 1_I)) \circ \alpha_{R,I,I}$ and the left and bottom is the summand $\rho \circ (f\ot1) \circ ((1\ot e)\ot1_I)$ we conclude that 
\begin{align*}
f \circ (1 \ot \rho) \circ (1\ot(e \ot 1_I)) \circ \alpha_{R,I,I} - \rho \circ (f\ot1) \circ ((1\ot e)\ot1_I)=0 
\end{align*}
So we are left with 
\begin{align*}
&
d^2(f) \circ ((1 \ot e) \ot e) \circ \alpha_{R,I,I}^{-1} \circ ((1,0) \ot \rho_I^{-1})
\\&=
(\mu \circ [1 \ot (f \circ (e \ot e))] \circ \alpha_{R,I,I} - f \circ [\rho \ot e]) \circ \alpha_{R,I,I}^{-1} \circ ((1,0) \ot \rho_I^{-1})
\\&=
 \mu \circ [(1,0) \ot (f \circ (e \ot e) \circ \rho_I^{-1})] - f \circ ((1,0) \ot e ) = 0 
\end{align*}
where we in the last summand have that $(\rho \ot 1)\alpha_{R,I,I}^{-1}(1 \ot \rho_I^{-1})=1_{R\ot I}$ by the coherence theorem. Hence in the verification of the right unit law we summarise and are left with 
\begin{align*}
&
\mu_{R\ltimes_f R} \circ (1_{R\ltimes_f R} \ot e_{R\ltimes_f R} ) \circ \rho_{R\ltimes_fR}^{-1}
\\&=
\begin{pmatrix}
\mu \circ ((1,0)\ot e) \\
f   \circ ((1,0)\ot e) - \mu \circ ((1,0) \ot (f \circ (e \ot e) \circ \rho_{I}^{-1}) + \mu \circ ((0,1)\ot e)
\end{pmatrix}\rho_{R\ltimes_f R}^{-1}
\\&=
\begin{pmatrix}
\mu \circ ((1,0)\ot e) \\
\mu \circ ((0,1)\ot e)
\end{pmatrix}\rho_{R\ltimes_f R}^{-1}
=
\begin{pmatrix}
\rho \circ ((1,0)\ot 1_I) \\
\rho \circ ((0,1)\ot 1_I)
\end{pmatrix}\rho_{R\ltimes_f R}^{-1}
=
\begin{pmatrix}
((1,0)\ot 1_I) \circ \rho_{R\ltimes_f R} \\
((0,1)\ot 1_I) \circ \rho_{R\ltimes_f R}
\end{pmatrix}\rho_{R\ltimes_f R}^{-1}
\\&=
\begin{pmatrix}
1&0\\
0&1  
\end{pmatrix}
=1_{R\ltimes_f R}
\end{align*}
where we have used the unitary law $\mu \circ ((1,0)\ot e) = \rho\circ ((1,0)\ot 1_I)$ and $\mu\circ((0,1)\ot e) = \rho \circ ((0,1)\ot 1_I)$, and that $\rho$ is a natural transformation. This proves the right unitary law. \mbox{}\hfill\qed\end{mypf}

\begin{myprop}
Let $(\fC,\ot,I,\alpha,\lambda,\rho)$ be an additive monoidal category and let $(R,\mu,e)$ be a ring object in $\fC$. If $\overline{f}=\overline{g}$ in $\HH^2(R)$ then $R\ltimes_f R \cong R\ltimes_g R$ as ring objects. \end{myprop}

\begin{mypf}
We have to construct an isomorphism of ring objects (see Definition \ref{def:mro}). Since $\overline{f}=\overline{g}$ there is a $h\in\Hom_{\fC}(R , R)$ such that $f=g+d^1(h)=g+(\mu \circ (1 \ot h) - h \circ \mu + \mu \circ (h \ot 1))$. We consider the following morphisms in $\fC$
\begin{align}
\phi&= 
\begin{pmatrix}
 1 & 0 \\
 h & 1 
\end{pmatrix}
: R \ltimes_f R \to R \ltimes_g R \\
\psi&= 
\begin{pmatrix}
 1 & 0 \\
-h & 1 
\end{pmatrix}
: R \ltimes_g R \to R \ltimes_f R. 
\end{align}
These morphisms are clearly mutual inverses of one another, 
\begin{align*}
\psi\phi= 
\begin{pmatrix}
 1 & 0 \\
-h & 1 
\end{pmatrix}
\begin{pmatrix}
 1 & 0 \\
 h & 1 
\end{pmatrix}
=
\begin{pmatrix}
 1   & 0 \\
-h+h & 1 
\end{pmatrix}
= 
\begin{pmatrix}
 1 & 0 \\
 0 & 1 
\end{pmatrix}
=1_{R\ltimes_f R}
\\
\phi\psi= 
\begin{pmatrix}
 1 & 0 \\
 h & 1 
\end{pmatrix}
\begin{pmatrix}
 1 & 0 \\
-h & 1 
\end{pmatrix}
=
\begin{pmatrix}
 1   & 0 \\
 h-h & 1 
\end{pmatrix}
= 
\begin{pmatrix}
 1 & 0 \\
 0 & 1 
\end{pmatrix}
=1_{R\ltimes_g R}
\end{align*}
So we have to check if $\phi$ and $\psi$ really define morphisms of ring objects. Starting with $\phi$ we have to show that the following two diagrams commute 
\newline
\begin{minipage}{.5\textwidth}
\begin{center}
\begin{tikzpicture}
\matrix(m)[matrix of math nodes,row sep=2.6em,column sep=2.8em,text height=1.5ex,text depth=0.25ex]
{
    (R\ltimes_f R)\ot(R\ltimes_f R)   &   (R\ltimes_g R)\ot(R\ltimes_g R) \\
     R\ltimes_f R                     &    R\ltimes_g R                   \\
};
\draw[ ->,font=\scriptsize](m-1-1) edge         node[above]{$\phi\ot\phi               $} (m-1-2);
\draw[ ->,font=\scriptsize](m-2-1) edge         node[above]{$\phi                      $} (m-2-2);
\draw[ ->,font=\scriptsize](m-1-1) edge         node[left ]{$\mu_{R\ltimes_f R}        $} (m-2-1);
\draw[ ->,font=\scriptsize](m-1-2) edge         node[right]{$\mu_{R\ltimes_g R}        $} (m-2-2);
\end{tikzpicture}
\end{center}
\end{minipage}
\begin{minipage}{.5\textwidth}
\begin{center}
\begin{tikzpicture}
\matrix(m)[matrix of math nodes,row sep=2.6em,column sep=2.8em,text height=1.5ex,text depth=0.25ex]
{
     I                                &        R\ltimes_f R               \\
     I                                &        R\ltimes_g R               \\
};
\draw[ ->,font=\scriptsize](m-1-1) edge         node[above]{$e_{R\ltimes_f R}          $} (m-1-2);
\draw[ ->,font=\scriptsize](m-2-1) edge         node[above]{$e_{R\ltimes_g R}          $} (m-2-2);
\draw[ ->,font=\scriptsize](m-1-1) edge         node[left ]{$1_I                       $} (m-2-1);
\draw[ ->,font=\scriptsize](m-1-2) edge         node[right]{$\phi                      $} (m-2-2);
\end{tikzpicture}
\end{center}
\end{minipage}
Starting with the left diagram, we have 
\begin{align*}
\mu_{R\ltimes_g R} \circ (\phi\ot\phi)
=&
\begin{pmatrix}
\mu &  0  &  0  & 0 \\
 g  & \mu & \mu & 0 
\end{pmatrix}
\begin{pmatrix}
(1,0)\ot(1,0)\\
(1,0)\ot(0,1)\\
(0,1)\ot(1,0)\\
(0,1)\ot(0,1)\\
\end{pmatrix}
\begin{pmatrix}
\begin{pmatrix}
1 & 0 \\
h & 1 
\end{pmatrix}
\ot
\begin{pmatrix}
1 & 0 \\
h & 1 
\end{pmatrix}
\end{pmatrix}
\\&=
\begin{pmatrix}
\mu \circ ((1,0)\ot(1,0))\\
g   \circ ((1,0)\ot(1,0)) + \mu \circ ((h+1,1)\ot(h+1,1))\\
\end{pmatrix}
\end{align*}
and 
\begin{align*}
\phi\circ\mu_{R\ltimes_f R}&=
\begin{pmatrix}
1 & 0 \\
h & 1 
\end{pmatrix}
\begin{pmatrix}
\mu &  0  &  0  & 0 \\
f   & \mu & \mu & 0 
\end{pmatrix}
\begin{pmatrix}
(1,0)\ot(1,0)\\
(1,0)\ot(0,1)\\
(0,1)\ot(1,0)\\
(0,1)\ot(0,1)\\
\end{pmatrix}
\\&=
\begin{pmatrix}
\mu \circ ((1,0)\ot(1,0))\\
h \circ \mu \circ ((1,0)\ot(1,0))+f \circ ((1,0)\ot(1,0))+\mu \circ ((1,1)\ot(1,1))\\
\end{pmatrix}
\\&=
\begin{pmatrix}
\mu \circ ((1,0)\ot(1,0))\\
g   \circ ((1,0)\ot(1,0)) + \mu \circ ((h+1,1)\ot(h+1,1))\\
\end{pmatrix}
\end{align*}
where the last equality is determined by the expression of $f$ with $g$ and $h$. We observe that the expressions are the same, and conclude that the left diagram commutes, and the multiplicative structure is preserved. For the right diagram, or the preservation of the unit
\begin{align*}
\phi \circ e_{R\ltimes_fR}=
\begin{pmatrix}
1 & 0 \\ 
h & 1
\end{pmatrix}
\begin{pmatrix}
e \\ 
-f \circ (e \ot e) \circ \lambda_I^{-1}
\end{pmatrix}
=
\begin{pmatrix}
e \\ 
-g \circ (e \ot e) \circ \lambda_I^{-1}
\end{pmatrix}=e_{R\ltimes_gR}
\end{align*}
where the middle equality is obtained from 
\begin{align*}
f \circ (e \ot e) &= g \circ (e \ot e)+(\mu \circ (e \ot [h \circ e]) - h \circ \mu \circ (e \ot e) + \mu \circ ( [h \circ e] \ot e ) ) \\&= g \circ (e \ot e) + h \circ e \circ \lambda_I^{-1}. 
\end{align*}
We conclude that $\phi$ is a morphism of ring objects. A similar justification proves that $\psi$ is a morphism of ring objects as well. \mbox{}\hfill\qed\end{mypf}

\section{The Hochschild cohomology ring} 
In this section we prove that $\HH^*(R)=\bigoplus_{i=0}^{\infty} \HH^i(R)$ is a graded-commutative ring with the cup product. In \cite{ger-63} and \cite{gs-86} this result is proved when $A$ is an associative ring, i.e.\ for classical Hochschild cohomology. Our proof for ring objects in $\Ab$-enriched monoidal categories is motivated by that of Gerstenhaber. From \cite{ger-63} and \cite{gs-86} we will recall some of the needed concepts such as right pre-Lie system and right pre-Lie ring, and make use of these in the proof. Throughout this section, let $(\fC,\ot,I,\alpha,\lambda,\rho)$ be an $\Ab$-enriched monoidal category and let $(R,\mu,e)$ be a ring object in $\fC$. 

\subsection{Graded rings and pre-Lie systems}
Recall that a ring $\Lambda$ is \emph{$\mathbb{Z}$-graded} if it decomposes as $\Lambda=\bigoplus_{i=0}^{\infty}\Lambda^i$ such that $\Lambda^i\Lambda^j\subseteq\Lambda^{i+j}$ for all $i,j \in \mathbb{Z}$. The elements of $\Lambda^i$ are said to be \emph{homogeneous of degree $i$}. For homogeneous elements $a\in\Lambda^{i}$ and $b\in\Lambda^{j}$ the \emph{graded commutator} is defined as $[a,b]=ab-(-1)^{ij}ba$. The ring $\Lambda$ is said to be \emph{graded-commutative} if $[a,b]=0$, i.e.\ $ab=(-1)^{ij}ba$, for all homogeneous elements $a \in \Lambda^i$ and $b \in \Lambda^j$. We say that $\Lambda$ is a \emph{graded right pre-Lie ring} if for homogeneous elements $a\in\Lambda^{i+1}$, $b\in\Lambda^{j+1}$ and $c\in\Lambda^{k+1}$ we have 
\begin{align*}
(ab)c-(-1)^{jk}(ac)b=a(bc-(-1)^{jk}cb). 
\end{align*}
By a \emph{right pre-Lie system} $\{V_m,\bullet_i\}$ we mean a system of objects $V_{m}$ together with an operation $\bullet_i={\bullet_i}_{(m,n)}:V_m\ot V_n \to V_{m+n-1}$ for $i\in\{1,2,\dots,m\}$ that for $f\in V_m$, $g\in V_n$ and $h\in V_p$ satisfies 
\begin{align}\label{al:preLie}
(f\bullet_i g)\bullet_j h = 
\left\{
\begin{aligned}
&(f \bullet_j  h) \bullet_{i+p-1} g     &\text{if} & &1 \leq j \leq i-1& \\
& f \bullet_i (g  \bullet_{j-i+1} h)    &\text{if} & &i \leq j \leq n  &   
\end{aligned}
\right.
\end{align}
Further, the following identifications can be useful. From the first property we have  
\begin{align}
(f \bullet_j h ) \bullet_{i+p-1} g =(f \bullet_{i+p-1} g ) \bullet_{j+n-1} h \qquad \text{ if } \qquad 1 \leq i+p-1 \leq j-1
\end{align}
Now reading the expression above from right to left replacing $i+p-1$ by $i$ and $j+n-1$ by $j$ we get (for $i \leq j-n+1-1 \leq m-1$)  
\begin{align}
(f\bullet_i g)\bullet_j h = ( f \bullet_{j-n+1} h ) \bullet_{i} g \qquad \text{ if } \qquad i + n  \leq j \leq m + n - 1 
\end{align}

For a right pre-Lie system $\{V_m,\bullet_i\}$ we define the \emph{composition product} of $f\in V_m$ and $g\in V_n$ by
\begin{align}
f \bullet g = \sum_{i=1}^{m}(-1)^{(i-1)(n-1)}f \bullet_i g.  
\end{align}
We observe that the composition product also gives an operation $\bullet=\bullet_{(m,n)}:V_m \ot V_n \to V_{m+n-1}$. 

We recall the following results for composition products (this result and its proof corresponds to \cite[Theorem 2 and Corollary]{ger-63}). 

\begin{myprop}\label{prop:preLie}
Given an arbitrary right pre-Lie system $\{V_m,\bullet_i\}$, and $f\in V_m$, $g\in V_n$ and $h\in V_p$, then 
\begin{itemize}
\item[(i)  ] \[(f \bullet g) \bullet h - f \bullet ( g \bullet h ) = \sum_{(1 \leq j \leq i-1 )\vee( n + i \leq j \leq m + n - 1)} (-1)^{(n-1)(i-1)+(p-1)(j-1)}(f \bullet_i g) \bullet_j h\]
\item[(ii) ] \[(f \bullet g ) \bullet h - f \bullet (g \bullet h) = (-1)^{(n-1)(p-1)}\left((f \bullet h) \bullet g - f \bullet (h \bullet g) \right)\]
\item[(iii)] Let $A=\amalg_m V_m$ be the coproduct of the objects $V_m$ extending the composite operation defined on homogeneous elements to an operation $\bullet : A \ot A \to A$, then $A$ becomes a right pre-Lie ring. 
\end{itemize}
\end{myprop}

\begin{mypf}
We adapt the proof from \cite{ger-63}, but skip some details when they occur in \cite{ger-63}. For (i) we write out the expressions on the left side, 
\begin{align*}
(f \bullet  g) \bullet h  = \sum^{m+n-1}_{j=1}\sum^{m}_{i=1}   (-1)^{(j  -1)(p  -1)+(i     -1)(n-1)}(f \bullet_i      g) \bullet_j        h \\
 f \bullet (g  \bullet h) = \sum^{m}_{\xi=1}\sum^{n}_{\omega=1}(-1)^{(\xi-1)(n+p-2)+(\omega-1)(p-1)} f \bullet_{\xi} (g  \bullet_{\omega} h) 
\end{align*}
By the second defining property of right pre-Lie systems we identify $(f\bullet_i g)\bullet_j h =  f \bullet_i (g  \bullet_{j-i+1} h )$ when $i \leq j \leq n $, so $(f\bullet_i g)\bullet_j h -  f \bullet_{\xi} (g  \bullet_{\omega} h )=0$ if $\xi=i$ and $\omega=j-1+1$. First we check if the signs from the sums above match. The expression $(f\bullet_i g)\bullet_j h$ occurs with the sign $(-1)^{(j-1)(p-1)+(i-1)(n-1)}$ while $f \bullet_i (g  \bullet_{j-i+1} h )$ occurs with the same sign 
\[(-1)^{(i-1)(n+p-2)+(j-i)(p-1)}=(-1)^{(j-1)(p-1)+(i-1)(n-1)}  \] 
hence the suggested terms in the sums cancel. Moreover every term in the second sum cancels against some term in the first sum. To see this we pick an arbitrary term $f \bullet_{\xi} (g  \bullet_{\omega} h)$ where $1 \leq \xi \leq m$ and $1 \leq \omega \leq n$. Again this arbitrary term cancels against $(f\bullet_i g)\bullet_j h$ whenever $\xi=i$ and $\omega=j-1+1$, i.e.\ terms where $1 \leq i \leq m$ and $i \leq j \leq n+i-1$. These terms are clearly in the first sum. Hence we are left with terms in the first sum with either $1 \leq j \leq i - 1 $ or $ n + i \leq j \leq m + n - 1$, and the assertion follows. 

For (ii), we have from (i) that 
\begin{align*}
(f \bullet g) \bullet h - f \bullet ( g \bullet h ) &=& \sum_{(1 \leq j      \leq i  -1 )\vee( n + i   \leq j      \leq m + n - 1)} (-1)^{(n-1)(i  -1)+(p-1)(j     -1)}(f \bullet_i     g) \bullet_j h \\
(f \bullet h) \bullet g - f \bullet ( h \bullet g ) &=& \sum_{(1 \leq \omega \leq \xi-1 )\vee( p + \xi \leq \omega \leq m + p - 1)} (-1)^{(p-1)(\xi-1)+(n-1)(\omega-1)}(f \bullet_{\xi} h) \bullet_{\omega} g
\end{align*}
First we examine the case where $1 \leq j\leq i-1$ in the first sum, and see if we can transform the terms $(f \bullet_i g) \bullet_j h$ to terms of the form $(f \bullet_{\xi} h) \bullet_{\omega} g$ in the second sum. In the defining property of right pre-Lie systems we have for $1 \leq j \leq i-1$ that $(f\bullet_i g)\bullet_j h=(f \bullet_j  h) \bullet_{i+p-1} g$, i.e.\ $(f\bullet_i g)\bullet_j h=(f \bullet_{\xi}  h) \bullet_{\omega} g$ when $\xi=j$ and $\omega=i+p-1$. The term $(f\bullet_i g)\bullet_j h$ occurs in the first sum with the sign $(-1)^{(j-1)(p-1)+(i-1)(n-1)}$ while $(f \bullet_j  h) \bullet_{i+p-1} g$ occurs in the second sum with the sign 
\[(-1)^{(p-1)(j-1)+(n-1)(i+p)}=(-1)^{(j-1)(p-1)+(i-1)(n-1)}(-1)^{(n-1)(p-1)}.\] 
Hence, as asserted, we have to multiply the suggested term in the second sum with $(-1)^{(n-1)(p-1)}$ in order to transform it to the corresponding term in the first sum. All required terms of the form $(f \bullet_j  h) \bullet_{i+p-1} g$ occur in the second sum, by $p + \xi \leq \omega \leq m + p - 1$ it follows that $1 \leq j \leq i-1 $. For the case where $ i + n \leq j \leq m + n - 1 $ we have by the defining property of right pre-Lie systems that $ ( f \bullet_i g ) \bullet_j h = ( f \bullet_{j-n+1} h ) \bullet_i g $, so we want to transform the term $ ( f \bullet_i g ) \bullet_j h$ in the first sum to a term of the form $(f \bullet_{\xi} h) \bullet_{\omega} g$ in the second sum with $\xi=j-n+1$ and $\omega=i$. This term in the first sum occurs with the sign $(-1)^{(n-1)(i-1)+(p-1)(j-1)}$ while the term in the second occurs with sign 
\begin{align*}
(-1)^{(p-1)(j-n)+(n-1)(i-1)}&=(-1)^{(n-1)(i-1)+(p-1)(j-1)}(-1)^{-(p-1)(n-1)}\\
                            &=(-1)^{(n-1)(i-1)+(p-1)(j-1)}(-1)^{(p-1)(n-1)}. 
\end{align*}
So the term $( f \bullet_i g ) \bullet_j h$ in the first sum can be transformed to the term $( f \bullet_{j-n+1} h) \bullet_i g $ in the second sum by multiplying it with $(-1)^{(p-1)(n-1)}$. These terms clearly occur in the second sum so from $1 \leq \omega \leq \xi -1$ we get $i+n \leq j \leq m+n-1$. With these identifications we get that $(f \bullet g) \bullet h - f \bullet ( g \bullet h ) = (-1)^{(n-1)(p-1)}((f \bullet h) \bullet g - f \bullet ( h \bullet g ))$ and the assertion follows. 

Finally, (iii) is a direct consequence of (ii), which can be rewritten as 
\begin{align*}
(f \bullet g ) \bullet h - (-1)^{(n-1)(p-1)}\left((f \bullet h) \bullet g \right) = f \bullet (g \bullet h) - (-1)^{(n-1)(p-1)}\left( f \bullet (h \bullet g) \right)
\end{align*}
and the defining property for a right pre-Lie ring follows. \mbox{}\hfill\qed\end{mypf}

\subsection{A pre-Lie system}
The next objective is to prove that the following construction involving ring objects is a right pre-Lie system. As in the definition of the Hochschild cochain complex we denote $C^m(R)=\Hom_{\fC}(R^{\ot m},R)$ when $m\geq1$. For $f\in C^m(R)$ and $g\in C^n(R)$ (both $m,n\geq1$) we define the operation $\bullet_i$ as the composition of 
\begin{align*}
R^{\ot(m+n-1)}\xrightarrow{\alpha^{i-1,n}_{m+n-1}}(R^{\ot(i-1)}\ot R^{\ot n})\ot R^{\ot(m-i)}\xrightarrow{(1_{R^{\ot(i-1)}}\ot g) \ot 1_{R^{(m-i)}}}(R^{\ot(i-1)}\ot R)\ot R^{\ot(m-i)}\\
\xrightarrow{(\alpha^{i-1,1}_m)^{-1}}R^{\ot m}\xrightarrow{f}R
\end{align*}
i.e.\ $g$ evaluated at the ``$i$th'' place (of $R^{\ot (m+n-1)}$) before composing with $f$. We write out this as  
\begin{align}\label{eq:monpreLie}
f \bullet_i g = f \circ \alpha^{i-1,1}_m \circ [(1_{R^{\ot(i-1)}}\ot g) \ot 1_{R^{(m-i)}}] \circ \alpha^{i-1,n}_{m+n-1}:R^{\ot(m+n-1)}\to R. 
\end{align}
Clearly $\bullet_i$ defines an operation $C^m(R) \otimes_{\mathbb{Z}} C^n(R) \to C^{m+n-1}(R)$. 

\begin{myprop}
The system  $\{C^m(R),\bullet_i\}$ is a pre-Lie system. \end{myprop}

\begin{mypf}
We have to check if the relations from Equation (\ref{al:preLie}) hold. Let $f\in C^m(R)$, $g\in C^n(R)$ and $h \in C^p(R)$ (with all $m,n,p\geq1$). For $1 \leq j \leq i-1$ we have to show that $(f\bullet_i g)\bullet_j h = (f \bullet_j  h) \bullet_{i+p-1}$. Consider the following diagram 
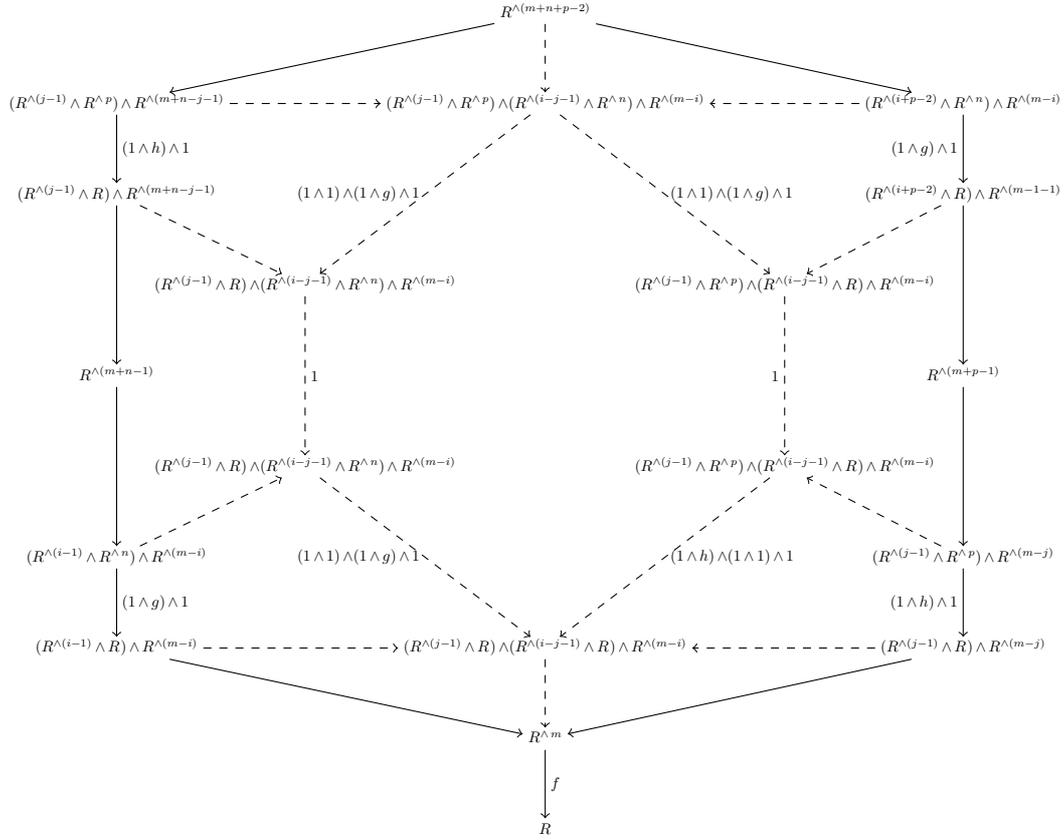
\begin{center}
\begin{tikzpicture}[scale=0.58, every node/.style={scale=0.58}, every ->/.style={scale=0.58}]
\matrix(m)[matrix of math nodes,row sep=2.6em,column sep=-3em,text height=1.5ex,text depth=0.25ex]
{   
                                                && R^{\ot(m+n+p-2)}                                                            &&                                        \\
 (R^{\ot(j-1)}\ot R^{\ot p})\ot R^{\ot(m+n-j-1)}&&(R^{\ot(j-1)}\ot R^{\ot p})\ot (R^{\ot(i-j-1)}\ot R^{\ot n}) \ot R^{\ot(m-i)}&                                                                          &(R^{\ot(i+p-2)}\ot R^{\ot n})\ot R^{\ot(m-i)}    \\
 (R^{\ot(j-1)}\ot R^{     })\ot R^{\ot(m+n-j-1)}&                                                                              &&                                                                         &(R^{\ot(i+p-2)}\ot R^{}     )\ot R^{\ot(m-1-1)}  \\
                                                & (R^{\ot(j-1)}\ot R^{     })\ot (R^{\ot(i-j-1)}\ot R^{\ot n}) \ot R^{\ot(m-i)}&&(R^{\ot(j-1)}\ot R^{\ot p})\ot (R^{\ot(i-j-1)}\ot R^{}) \ot R^{\ot(m-i)} &                                                 \\
  R^{\ot(m+n-1)}                                &                                                                              &&                                                                         & R^{\ot(m+p-1)}                                  \\
                                                & (R^{\ot(j-1)}\ot R^{     })\ot (R^{\ot(i-j-1)}\ot R^{\ot n}) \ot R^{\ot(m-i)}&&(R^{\ot(j-1)}\ot R^{\ot p})\ot (R^{\ot(i-j-1)}\ot R^{}) \ot R^{\ot(m-i)} &                                                 \\
 (R^{\ot(i-1)}\ot R^{\ot n})\ot R^{\ot(m-i)}    &                                                                              &&                                                                         &(R^{\ot(j-1)}\ot R^{\ot p})\ot R^{\ot(m-j)}      \\
 (R^{\ot(i-1)}\ot R^{     })\ot R^{\ot(m-i)}    &&(R^{\ot(j-1)}\ot R^{     })\ot (R^{\ot(i-j-1)}\ot R^{     }) \ot R^{\ot(m-i)}&&(R^{\ot(j-1)}\ot R^{     })\ot R^{\ot(m-j)}     \\
                                                && R^{\ot m       }                                                            &&                                        \\
                                                && R^{            }                                                            &&                                        \\
};
\draw[ ->](m-1-3) edge         node[left ]{$                          $} (m-2-1);
\draw[ ->](m-1-3) edge[dashed] node[left ]{$                          $} (m-2-3);
\draw[ ->](m-1-3) edge         node[left ]{$                          $} (m-2-5);
\draw[ ->](m-2-1) edge[dashed] node[left ]{$                          $} (m-2-3);
\draw[ ->](m-2-5) edge[dashed] node[left ]{$                          $} (m-2-3);
\draw[ ->](m-2-1) edge         node[right]{$(1 \ot h ) \ot 1          $} (m-3-1);
\draw[ ->](m-2-3) edge[dashed] node[left ]{$(1 \ot 1)\ot(1\ot g) \ot 1$} (m-4-2);
\draw[ ->](m-2-3) edge[dashed] node[right]{$(1 \ot 1)\ot(1\ot g) \ot 1$} (m-4-4);
\draw[ ->](m-2-5) edge         node[left ]{$(1 \ot g ) \ot 1          $} (m-3-5);
\draw[ ->](m-3-1) edge[dashed] node[left ]{$                          $} (m-4-2);
\draw[ ->](m-3-5) edge[dashed] node[left ]{$                          $} (m-4-4);
\draw[ ->](m-3-1) edge         node[right]{$                          $} (m-5-1);
\draw[ ->](m-5-1) edge         node[right]{$                          $} (m-7-1);
\draw[ ->](m-4-2) edge[dashed] node[right]{$1                         $} (m-6-2);
\draw[ ->](m-4-4) edge[dashed] node[left ]{$1                         $} (m-6-4);
\draw[ ->](m-3-5) edge         node[left ]{$                          $} (m-5-5);
\draw[ ->](m-5-5) edge         node[left ]{$                          $} (m-7-5);
\draw[ ->](m-7-1) edge[dashed] node[left ]{$                          $} (m-6-2);
\draw[ ->](m-7-5) edge[dashed] node[left ]{$                          $} (m-6-4);
\draw[ ->](m-7-1) edge         node[right]{$(1 \ot g ) \ot 1          $} (m-8-1);
\draw[ ->](m-6-2) edge[dashed] node[left ]{$(1 \ot 1)\ot(1\ot g) \ot 1$} (m-8-3);
\draw[ ->](m-6-4) edge[dashed] node[right]{$(1 \ot h)\ot(1\ot 1) \ot 1$} (m-8-3);
\draw[ ->](m-7-5) edge         node[left ]{$(1 \ot h ) \ot 1          $} (m-8-5);
\draw[ ->](m-8-1) edge[dashed] node[left ]{$                          $} (m-8-3);
\draw[ ->](m-8-5) edge[dashed] node[left ]{$                          $} (m-8-3);
\draw[ ->](m-8-1) edge         node[right]{$                          $} (m-9-3);
\draw[ ->](m-8-3) edge[dashed] node[right]{$                          $} (m-9-3);
\draw[ ->](m-8-5) edge         node[left ]{$                          $} (m-9-3);
\draw[ ->](m-9-3) edge         node[right]{$f                         $} (m-10-3);
\end{tikzpicture}
\end{center}
where the composition of the solid arrows on the left is $(f \bullet_i g)\bullet_jh$ while the composition of the solid arrows on the right is $(f\bullet_jh)\bullet_{i+p-1}g$. Unlabeled arrows are associativity relations. Dashed arrows are added for the purpose that it is easer to extract information. We also remark that $i-j-1\geq0$, so in the case of $i-j-1=0$ we have that $R^{\ot i-j-1}=R^{\ot 0}$ is the empty symbol. The middle hexagon commutes by successive compositions. The rest of the diagram commutes by naturality and coherence. Hence we conclude that $(f \bullet_i g)\bullet_jh=(f\bullet_jh)\bullet_{i+p-1}g$ when $1\leq j \leq i-1$. 

For the other defining property when $i\leq j \leq n$ we want to show that $ (f\bullet_i g)\bullet_j h = f \bullet_i (g  \bullet_{j-i+1} h) $. Consider the following diagram 
\begin{center}
\begin{tikzpicture}[scale=0.6, every node/.style={scale=0.6}, every ->/.style={scale=0.6}]
\matrix(m)[matrix of math nodes,row sep=2.6em,column sep=2.8em,text height=1.5ex,text depth=0.25ex]
{   
  R^{\ot(m+n+p-2)}                              && R^{\ot(m+n+p-2)}                                                                              \\
                                                &&(R^{\ot(i-1)}\ot   R^{\ot(n+p-1)})\ot R^{\ot(m-i)}                                              \\
 (R^{\ot(j-1)}\ot R^{\ot p})\ot R^{\ot(m+n-j-1)}&&(R^{\ot(i-1)}\ot ((R^{\ot(j-i  )} \ot R^{\ot p   })\ot R^{\ot(n+i-j-1)}))\ot R^{\ot(m-i)}       \\
 (R^{\ot(j-1)}\ot R^{     })\ot R^{\ot(m+n-j-1)}&&(R^{\ot(i-1)}\ot ((R^{\ot(j-i  )} \ot R^{        })\ot R^{\ot(n+i-j-1)}))\ot R^{\ot(m-i)}       \\
  R^{\ot(m+n-1)}                                &&                                                                                               \\
 (R^{\ot(i-1)}\ot R^{\ot n})\ot R^{\ot(m-i)}    &&(R^{\ot(i-1)}\ot R^{\ot n})\ot R^{\ot(m-i)}                                                    \\
 (R^{\ot(i-1)}\ot R^{     })\ot R^{\ot(m-i)}    &&(R^{\ot(i-1)}\ot R^{     })\ot R^{\ot(m-i)}                                                    \\
  R^{\ot m}                                     && R^{\ot m}                                                                                     \\
  R                                             && R                                                                                             \\
};
\draw[ ->](m-1-1) edge         node[right]{$                          $} (m-3-1);
\draw[ ->](m-3-1) edge         node[right]{$(1 \ot h)\ot 1            $} (m-4-1);
\draw[ ->](m-4-1) edge         node[right]{$                          $} (m-5-1);
\draw[ ->](m-5-1) edge         node[right]{$                          $} (m-6-1);
\draw[ ->](m-6-1) edge         node[right]{$(1 \ot g)\ot 1            $} (m-7-1);
\draw[ ->](m-7-1) edge         node[right]{$                          $} (m-8-1);
\draw[ ->](m-8-1) edge         node[right]{$f                         $} (m-9-1);
\draw[ ->](m-1-3) edge         node[left ]{$                          $} (m-2-3); 
\draw[ ->](m-2-3) edge         node[left ]{$                          $} (m-3-3); 
\draw[ ->](m-3-3) edge         node[left ]{$(1\ot((1\ot h)\ot 1))\ot1 $} (m-4-3); 
\draw[ ->](m-4-3) edge         node[left ]{$                          $} (m-6-3); 
\draw[ ->](m-6-3) edge         node[left ]{$(1 \ot g)\ot 1            $} (m-7-3); 
\draw[ ->](m-7-3) edge         node[left ]{$                          $} (m-8-3); 
\draw[ ->](m-8-3) edge         node[left ]{$f                         $} (m-9-3); 
\draw[ ->](m-1-1) edge[dashed] node[above]{$1                         $} (m-1-3);
\draw[ ->](m-3-1) edge[dashed] node[above]{$                          $} (m-3-3);
\draw[ ->](m-4-1) edge[dashed] node[above]{$                          $} (m-4-3);
\draw[ ->](m-6-1) edge[dashed] node[above]{$1                         $} (m-6-3);
\draw[ ->](m-7-1) edge[dashed] node[above]{$1                         $} (m-7-3);
\draw[ ->](m-8-1) edge[dashed] node[above]{$1                         $} (m-8-3);
\draw[ ->](m-9-1) edge[dashed] node[above]{$1                         $} (m-9-3);
\end{tikzpicture}
\end{center}
where the solid left side represents $(f \bullet_i g)\bullet_jh$ and the solid right side represents $f\bullet_i(g\bullet_{j-i+1}h)$. The unlabeled arrows are again associative relations. The dashed arrows are added for help reasons. Considering the horizontal dashed associative relations (unlabeled horizontal arrows) clearly $R^{\ot p}$ is placed in the same ``location'' in the tuple on both sides, so these associative relations contain some instance of $R^{\ot(j-1)}\to R^{\ot(i-1)}\ot R^{\ot(j-i)}$, which indeed makes perfect sense since $i\leq j$. Further the associative relations also contain some instance of $R^{\ot(m+n-j-1)}\to R^{\ot(n+i-j-1)}\ot R^{\ot(m-i)}$, which make sense since $j-1\leq n$. The top square of the diagram commutes by the coherence theorem, the second square from the top commutes by naturality, while the third commutes by the coherence theorem again. The rest of the diagram commutes clearly. Hence also the second relation is satisfied, and we conclude that the given construction is a right pre-Lie system. \mbox{}\hfill\qed\end{mypf}

\subsection{Graded commutativity of Hochschild cohomology} 
We now prove that $\HH^*(R)=\bigoplus_{i=0}^{\infty}\HH^i(R)$ is a graded-commutative ring. The multiplicative structure is given by the cup product which we now define. Let $f\in C^{m}(R)$ and $g\in C^n(R)$. We define the \emph{cup product} $f \smile g$ as the composition 
\begin{align}
f \smile g: 
\left\{
\begin{aligned}
&R^{\ot(m+n)}\xrightarrow{\alpha^{m,n}_{m+n}}       R^{\ot m}\ot R^{\ot n} \xrightarrow{f \ot g} R\ot R \xrightarrow{\mu} R& &\text{if}&  &m,n \geq 1& \\
&R^{\ot(m  )}\xrightarrow{\rho_{R^{\ot m}}^{-1}}    R^{\ot m}\ot I         \xrightarrow{f \ot g} R\ot R \xrightarrow{\mu} R& &\text{if}&  &m\geq1,n=0& \\
&R^{\ot(m+n)}\xrightarrow{\lambda_{R^{\ot n}}^{-1}} I        \ot R^{\ot n} \xrightarrow{f \ot g} R\ot R \xrightarrow{\mu} R& &\text{if}&  &m=0,n\geq1& \\
&I\xrightarrow{\lambda_I^{-1}=\rho_I^{-1}}          I \ot I                \xrightarrow{f \ot g} R\ot R \xrightarrow{\mu} R& &\text{if}&  &m=0,n=0   & \\
\end{aligned}
\right.
\end{align}
We remark that $\alpha^{m,n}_{m+n}=\alpha^{m,0}_{m+n}$, and this might happen to be the identity. We write the compositions together and get the following forms 
\begin{align}
f \smile g=
\left\{
\begin{aligned}
& \mu \circ (f \ot g) \circ \alpha^{m,n}_{m+n}       & &\text{if}& &m,n \geq 1& \\
& \mu \circ (f \ot g) \circ \rho_{R^{\ot m}}^{-1}    & &\text{if}& &m\geq1,n=0& \\
& \mu \circ (f \ot g) \circ \lambda_{R^{\ot n}}^{-1} & &\text{if}& &m=0,n\geq1& \\
& \mu \circ (f \ot g) \circ \lambda_I^{-1}           & &\text{if}& &m=0,n=0   &   
\end{aligned}
\right.
\end{align} 
Thus the cup product gives an operation $C^m(R) \otimes_{\mathbb{Z}} C^n(R) \to C^{m+n}(R). $

\begin{mylem}
The multiplicative unit $e:I \to R$ is a left and right unit for the cup product. \end{mylem}

\begin{mypf}
For the right unitary law consider the following diagram 
\begin{center}
\begin{tikzpicture}
\matrix(m)[matrix of math nodes,row sep=2.6em,column sep=5.8em,text height=1.5ex,text depth=0.25ex]
{
  R^{\ot m}    &  R^{\ot m}     \ot I &  R^{\ot m}     \ot  R             \\
  R            &  R             \ot I &  R             \ot  R             \\
               &                      &  R                                \\
};
\draw[ ->,font=\scriptsize](m-1-1) edge         node[above]{$\rho_{R^{\ot m}}^{-1}                    $} (m-1-2);
\draw[ ->,font=\scriptsize](m-1-2) edge         node[above]{$1_{R^{\ot m}}\ot e                       $} (m-1-3);
\draw[ ->,font=\scriptsize](m-1-1) edge         node[right]{$f                                        $} (m-2-1);
\draw[ ->,font=\scriptsize](m-1-2) edge         node[right]{$f \ot 1_I                                $} (m-2-2);
\draw[ ->,font=\scriptsize](m-1-3) edge         node[right]{$f \ot 1_R                                $} (m-2-3);
\draw[ ->,font=\scriptsize](m-2-1) edge         node[above]{$\rho_R^{-1}                              $} (m-2-2);
\draw[ ->,font=\scriptsize](m-2-2) edge         node[above]{$1_{R}\ot e                               $} (m-2-3);
\draw[ ->,font=\scriptsize](m-2-3) edge         node[right]{$\mu                                      $} (m-3-3);
\draw[ ->,font=\scriptsize](m-2-1) edge         node[below]{$1_{R}                                    $} (m-3-3);
\end{tikzpicture}
\end{center}
The bottom triangle commutes by the right unitary law, the left square commutes by naturality and the right by successive compositions. We conclude that 
\begin{align*}
f \smile e = \mu \circ (f \ot e) \circ \rho_{R^{\ot m}}^{-1} = \mu \circ (f \ot 1_R) \circ (1_{R^{\ot m}} \ot e) \circ \rho_{R^{\ot m}}^{-1} = f
\end{align*}
where the last equality follows from the outer part of the diagram. Similarly, $e$ is also a left unit for the cup product by the left unitary law. \mbox{}\hfill\qed\end{mypf}

We shall prove that $\HH^*(R)$ is graded-commutative, i.e.\ for $\overline{f}\in \HH^m(R)$ and $\overline{g}\in \HH^n(R)$ we show that $\overline{f} \smile \overline{g} = (-1)^{nm}\overline{g} \smile \overline{f}$. For this purpose let $f\in C^m(R)$ and $g \in C^n(R)$ with $m,n \geq 1$. Observe by a direct computation that the cup product can be given as 
\begin{align}\label{eq:cupcirc}
f\smile g = (\mu \bullet_1 f)\bullet_{m+1} g. 
\end{align}
where the operation $\bullet_i$ is as in the previous section. Next we claim that the Hochschild differential can be written as 
\begin{align}\label{eq:circdiff}
d^m(f)=-(f \bullet \mu - (-1)^{m-1} \mu \bullet f). 
\end{align}
To see this, note that 
\begin{align*}
- f \bullet \mu = - \sum^{m}_{i=1}(-1)^{(i-1)\cdot 1} f \bullet_i \mu = - \sum^{m}_{i=1} (-1)^{i-1} f \circ \alpha^{i-1,1}_m \circ \mu^{i-1}_{m+1} \circ \alpha^{i-1,2}_{m+1} 
\end{align*}
where we have used the identity $\mu^{i-1}_{m+1}=( 1_R^{\ot (i-1)}{}\ot \mu ) \ot 1_{R^{\ot(m-i)}}$. Furthermore,  
\begin{align*}
(-1)^{m-1} \mu \bullet f &= (-1)^{m-1} \sum_{i=1}^2(-1)^{(i-1)(m-1)}\mu \bullet_{i} f \\
                         &= (-1)^{m-1}(\mu \circ (f\ot1_R)+(-1)^{m-1}\mu \circ (1_R \ot f) \circ \alpha^{0,1}_{m+1})\\
                         &= (-1)^{m+1} \mu \circ (f\ot1_R)+          \mu \circ (1_R \ot f) \circ \alpha^{0,1}_{m+1}
\end{align*}
where we have used the identities $\alpha^{i-1,1}_2=1_{R \ot R}$ and $\alpha^{1,m}_{m+1}=\alpha^{0,1}_{m+1}$ (by the coherence theorem and notation). Written together this gives 
\begin{align*}
-(f \bullet \mu - (-1)^{m-1} \mu \bullet f)&=\mu \circ (1_R \ot f) \circ \alpha^{0,1}_{m+1}+\sum^{m}_{i=1} (-1)^{i} f \circ \alpha^{i-1,1}_m \circ \mu^{i-1}_{m+1} \circ \alpha^{i-1,2}_{m+1} \\&\mbox{ }+ (-1)^{m+1} \mu \circ (f\ot1_R) \\
                                       &= d^m(f)
\end{align*}
The following result corresponds to \cite[Theorem 3 and Corollary 1]{ger-63}).  

\begin{myprop}\label{prop:circdiff}
Let $(\fC,\ot,I,\alpha,\lambda,\rho)$ be an $\Ab$-enriched monoidal category and let $(R,\mu,e)$ be a ring object in $\fC$. Moreover let $f \in C^m(R) = \Hom_{\fC}(R^{\ot m},R)$ and $g \in C^n(R) = \Hom_{\fC}(R^{\ot n},R)$ where $m,n\geq1$. Then 
\begin{itemize}
\item[(i)  ] $f \bullet (d^ng) - d^{m+n-1}(f \bullet g) + (-1)^{n-1}(d^mf) \bullet g = (-1)^{n-1}\left(g \smile f - (-1)^{mn}f \smile g\right)$
\item[(ii) ] if $f$ and $g$ are cocycles (i.e.\ $d^mf=0$ and $d^ng=0$) then 
\[d^{m+n-1}(f \bullet g) = (-1)^{n}\left(g \smile f - (-1)^{mn}f \smile g\right). \]
\end{itemize}\end{myprop}

\begin{mypf}
For (i) we use the expression for the differential in Equation (\ref{eq:circdiff}) and write the three terms on the left as
\begin{align*}
                     f \bullet (d^ng)&= -              f \bullet ( g   \bullet \mu ) + (-1)^{n-1}    f   \bullet ( \mu \bullet g )           \\ 
 -         d^{m+n-1}(f \bullet g)    &=              ( f \bullet   g ) \bullet \mu   - (-1)^{m+n-2}  \mu \bullet ( f   \bullet g ))          \\ 
(-1)^{n-1} (d^mf)      \bullet g     &= - (-1)^{n-1} ( f \bullet \mu ) \bullet g     + (-1)^{m+n-2} (\mu \bullet   f ) \bullet g         
\end{align*}
By the defining property for right pre-Lie rings we have that 
\begin{align*}
( f \bullet g ) \bullet \mu - (-1)^{n-1} (f \bullet \mu ) \bullet g = f \bullet ( g \bullet \mu - (-1)^{n-1} \mu \bullet g)
\end{align*}
leaving us with 
\begin{align*}
(-1)^{m+n-2}((\mu \bullet f) \bullet g - \mu \bullet (f \bullet g))
\end{align*}
in the original expression. From Proposition \ref{prop:preLie} part (i) we have 
\begin{align*}
(\mu \bullet f) \bullet g - \mu \bullet ( f \bullet g ) = \sum_{(1 \leq j \leq i-1 )\vee( m + i \leq j \leq m + 1)} (-1)^{(m-1)(i-1)+(n-1)(j-1)}(\mu \bullet_i f) \bullet_j g
\end{align*}
where the only possible values for $i$ are $1$ (which gives $j=m+1$) and $i=2$ (which gives $j=1$), hence 
\begin{align*}
(\mu \bullet f) \bullet g - \mu \bullet ( f \bullet g ) = (-1)^{(n-1)m}(\mu \bullet_1 f) \bullet_{m+1} g + (-1)^{m-1} (\mu \bullet_2 f) \bullet_{1} g . 
\end{align*}
We identify (via the identification in Equation (\ref{eq:cupcirc}))
\begin{align*}
(\mu \bullet_1 f) \bullet_{m+1} g &= f \smile g \\
(\mu \bullet_2 f) \bullet_{1}   g &= (\mu \bullet_1 g) \bullet_{n+1} g = g \smile f 
\end{align*}
where the second equality follows from the defining property of right pre-Lie systems. Inserted back we get 
\begin{align*}
(-1)^{m+n-2}\left((-1)^{m-1}g \smile f-(-1)^{(n-1)m+1}f \smile g\right) = (-1)^{n-1}\left(g \smile f-(-1)^{nm}f \smile g\right)
\end{align*}
which proves part (i). 

For part (ii), when $d^mf=0$ and $d^ng=0$ we have $f \bullet (d^ng)=0$ and $(d^mf) \bullet g = 0$ by construction, hence 
\begin{align*}
d^{m+n-1}(f \bullet g) = (-1)^{n}\left(g \smile f - (-1)^{mn}f \smile g\right)
\end{align*}
as asserted. \mbox{}\hfill\qed\end{mypf}

We prove that the cohomology ring is graded-commutative. 

\begin{mythm}
Let $(\fC,\ot,I,\alpha,\lambda,\rho)$ be an $\Ab$-enriched category and let $(R,\mu,e)$ be a ring object in $\fC$. Then the Hochschild cohomology ring $\HH^*(R)=\bigoplus_{i=0}^{\infty}\HH^i(R)$ is graded-commutative, that is if $\overline{f}\in \HH^m(R)$ and $\overline{g}\in \HH^n(R)$ then $\overline{f} \smile \overline{g} = (-1)^{mn}(\overline{g} \smile \overline{f})$. Moreover, its unit in $\overline{e}$. \end{mythm}

\begin{mypf}
We start by examining the case where $m=0$ and $n\geq1$. Recall that if $\overline{f} \in \HH^0(R) = \Ker d^0$ then $f$ is in the centre of $R$ (as discussed in Section \ref{sec:centre}), and for any $\overline{g} \in \HH^n(R)$ consider the following diagram 
\begin{center}
\begin{tikzpicture}
\matrix(m)[matrix of math nodes,row sep=2.6em,column sep=2.6em,text height=1.5ex,text depth=0.25ex]
{   
                  &          & R^{\ot n} &         &                 \\
        R^{\ot n} &          & R         &         & R^{\ot n}       \\
  I \ot R^{\ot n} &  I \ot R &           & R \ot I & R^{\ot n} \ot I \\
  R \ot R^{\ot n} &  R \ot R &           & R \ot R & R^{\ot n} \ot R \\
                  &          & R         &         &                 \\
};
\draw[ ->,font=\scriptsize](m-1-3) edge         node[left ]{$1_{R^{\ot n}}             $} (m-2-1);
\draw[ ->,font=\scriptsize](m-1-3) edge         node[left ]{$g                         $} (m-2-3);
\draw[ ->,font=\scriptsize](m-1-3) edge         node[right]{$1_{R^{\ot n}}             $} (m-2-5);
\draw[ ->,font=\scriptsize](m-2-1) edge         node[left ]{$\lambda^{-1}_{R^{\ot n}}  $} (m-3-1);
\draw[ ->,font=\scriptsize](m-2-3) edge         node[left ]{$\lambda^{-1}_R            $} (m-3-2);
\draw[ ->,font=\scriptsize](m-2-3) edge         node[right]{$\rho^{-1}_R               $} (m-3-4);
\draw[ ->,font=\scriptsize](m-2-5) edge         node[right]{$\rho^{-1}_{R^{\ot n}}     $} (m-3-5);
\draw[ ->,font=\scriptsize](m-3-1) edge         node[left ]{$f\ot1_{R^{\ot n}}         $} (m-4-1);
\draw[ ->,font=\scriptsize](m-3-2) edge         node[left ]{$f\ot1_{R}                 $} (m-4-2);
\draw[ ->,font=\scriptsize](m-3-4) edge         node[right]{$1_R\ot f                  $} (m-4-4);
\draw[ ->,font=\scriptsize](m-3-5) edge         node[right]{$1_{R^{\ot n}}\ot f        $} (m-4-5);
\draw[ ->,font=\scriptsize](m-4-2) edge         node[above]{$\mu                       $} (m-5-3);
\draw[ ->,font=\scriptsize](m-4-4) edge         node[above]{$\mu                       $} (m-5-3);
\draw[ ->,font=\scriptsize](m-2-1) edge         node[above]{$g                         $} (m-2-3);
\draw[ ->,font=\scriptsize](m-2-5) edge         node[above]{$g                         $} (m-2-3);
\draw[ ->,font=\scriptsize](m-3-1) edge         node[above]{$1_I \ot g                 $} (m-3-2);
\draw[ ->,font=\scriptsize](m-3-5) edge         node[above]{$g   \ot 1_I               $} (m-3-4);
\draw[ ->,font=\scriptsize](m-4-1) edge         node[above]{$1_R \ot g                 $} (m-4-2);
\draw[ ->,font=\scriptsize](m-4-5) edge         node[above]{$g   \ot 1_R               $} (m-4-4);
\end{tikzpicture}
\end{center}
The diagram clearly commutes by the definition of the centre (the middle hexagon), naturality (squares on the left and on the right) and successive compositions. By composing on the left and the right side we have, respectively,  
\begin{align*}
\mu \circ (1_R \ot g   ) \circ (f   \ot 1_R) \circ \lambda^{-1}_{R^{\ot n}} = \mu \circ ( f \ot g ) \circ \lambda^{-1}_{R^{\ot n}} =& f \smile g \\
\mu \circ ( g  \ot 1_R ) \circ (1_R \ot f  ) \circ \rho^{-1}_{R^{\ot n}}    = \mu \circ ( g \ot f ) \circ \rho^{-1}_{R^{\ot n}}    =& g \smile f 
\end{align*}
and by the commutativity of the diagram we conclude that $f \smile g = (-1)^{0 \cdot n} g \smile f = g \smile f $. A similar argument holds when $\overline{f}\in \HH^{m}(R)$ ($m\geq1$) and $\overline{g} \in \HH^0(R)$, we replace $f$ by $g$ (and $R^{\ot n}$ by $R^{\ot m}$). When $\overline{f}\in \HH^{0}(R)$ and $\overline{g} \in \HH^0(R)$ replacing $R^{\ot n}$ with $I$ in the diagram above, and the same result follows. 

For the remaining case where $m,n \geq 1$ we apply Proposition \ref{prop:circdiff}. When $\overline{f}\in \HH^{m}(R)$ and $\overline{g} \in \HH^{n}(R)$ clearly every representative for $\overline{f}$ and $\overline{g}$ are cocycles (in $\Ker d^m$ and $\Ker d^n$, respectively) by definition. Hence we can apply part (ii) of the proposition, which states that 
\begin{align*}
d^{m+n-1}(f \bullet g) = (-1)^{n}\left(g \smile f - (-1)^{mn}f \smile g\right). 
\end{align*}
Moreover recall that $f \bullet g \in C^{m+n-1}(R)$, hence clearly $d^{m+n-1}(f \bullet g)$ is a coboundary (i.e.\ $d^{m+n-1}(f \bullet g) \in \IM d^{m+n-1}\subseteq C^{m+n}(R)$) and so vanishes in $\HH^{m+n}(R)$, that is $\overline{d^{m+n-1}(f \bullet g)}=\overline{0}$. For $\overline{f}\in \HH^m(R)$ and $\overline{g}\in \HH^n(R)$ we have then that 
\begin{align*}
\overline{0} = (-1)^{n}\left(\overline{g} \smile \overline{f} - (-1)^{mn}\overline{f} \smile \overline{g}\right), 
\end{align*}
hence $\overline{g} \smile \overline{f} = (-1)^{mn}\overline{f} \smile \overline{g}$ and the result follows. \mbox{}\hfill\qed\end{mypf}

In \cite{ger-63} Gerstenhaber prove that $\HH^*(A)$ is a \emph{Gerstenhaber algebra} with the \emph{bracket multiplication} $[f,g]=f \bullet g - (-1)^{nm} g \bullet f$ (where $f\in C^m(R)$ and $g\in C^n(R)$). Since the proof that of Gerstenhaber only make use of the pre-Lie system, we also conclude that $\HH^*(R)$ is a Gerstenhaber algebra. 

\mbox{}\\
\noindent\textbf{Acknowledgements.} I thank my advisor Petter A. Bergh and colleagues at the Department of Mathematical Sciences (NTNU) for valuable discussions throughout this project. The Research Council of Norway support my phd-project through the project \emph{Triangulated categories in algebra} (NFR 221893). 
 
\addcontentsline{toc}{section}{References} 
\printbibliography

{\indent\textsc{Magnus Hellstrøm-Finnsen, Institutt for matematiske fag, NTNU, N-7491 Trondheim, Norway}}\\
{\indent{\mbox{}} Email address: \url{magnuhel@math.ntnu.no}}

\end{document}